\documentclass[a4paper,10pt]{amsart}

\pdfoutput=1

\usepackage[T1]{fontenc}
\usepackage[utf8]{inputenc}
\usepackage{newunicodechar}

\usepackage[usenames,dvipsnames,svgnames,x11names,table]{xcolor}

\usepackage[]{accents,
amsmath,
amsfonts,
amssymb,
amsthm,
amscd,
amsxtra,
amstext,
latexsym,
syntonly,
tabularx}

\usepackage{pgf,tikz}
\usetikzlibrary{arrows}
\usepackage{enumitem}
\usepackage{graphicx}
\usepackage{mathrsfs}
\usepackage{stmaryrd}
\usepackage[normalem]{ulem}
\usepackage[all]{xy}

\usepackage[pagebackref, colorlinks=true, linkcolor=
Firebrick4,
citecolor=
Turquoise4,
filecolor=SteelBlue4,
urlcolor=
SteelBlue4
]{hyperref}

\theoremstyle{plain}
\newtheorem{xx}{xx}[section]
\newtheorem{thm}[xx]{Theorem}%[section]
\newtheorem*{thm*}{Theorem}
\newtheorem*{thmA*}{Theorem A}
\newtheorem*{thmB*}{Theorem B}
\newtheorem*{thmC*}{Theorem C}
\newtheorem{prop}[xx]{Proposition}%[theorem]
\newtheorem{cor}[xx]{Corollary}%[theorem]
\newtheorem*{cor*}{Corollary}%[theorem]
\newtheorem*{corA*}{Corollary A}%[theorem]
\newtheorem*{corB*}{Corollary B}%[theorem]
\newtheorem{lem}[xx]{Lemma}%[theorem]
\newtheorem*{lem*}{Lemma}%[theorem]

\theoremstyle{definition}
\newtheorem{defn}[xx]{Definition}%[theorem]
\newtheorem*{defn*}{Definition}%[theorem]
\newtheorem{ex}[xx]{Example}%[theorem]
\theoremstyle{remark}
\newtheorem{rem}[xx]{Remark}
\newtheorem*{rem*}{Remark}

\numberwithin{equation}{xx}

\DeclareMathOperator{\Add}{Add}

\DeclareMathOperator{\Arr}{Arr}

\DeclareMathOperator{\Char}{char}

\DeclareMathOperator{\coker}{coker}
\DeclareMathOperator{\Def}{Def}

\DeclareMathOperator{\depth}{depth}
\DeclareMathOperator{\Der}{Der}
\DeclareMathOperator{\cE}{E}

\DeclareMathOperator{\End}{End}
\DeclareMathOperator{\Ext}{Ext}
\DeclareMathOperator{\uExt}{\ul{\Ext}}

\DeclareMathOperator{\grade}{grade}

\DeclareMathOperator{\cH}{H}

\DeclareMathOperator{\Hom}{Hom}
\DeclareMathOperator{\uHom}{\ul{\Hom}}
\DeclareMathOperator{\sHom}{\ms{H}\hspace{-0.15em}\mathnormal{om}}

\DeclareMathOperator{\id}{{id}}
\DeclareMathOperator{\im}{im}
\DeclareMathOperator{\Injdim}{inj{.}dim}

\DeclareMathOperator{\Mod}{mod}

\DeclareMathOperator{\ob}{ob}
\DeclareMathOperator{\pdim}{pdim}

\DeclareMathOperator{\Quot}{Quot}

\DeclareMathOperator{\rk}{rk}
\DeclareMathOperator{\Sets}{\cat{Sets}}

\DeclareMathOperator{\Spec}{Spec}

\DeclareMathOperator{\Sym}{Sym}
\DeclareMathOperator{\Syz}{Syz}

\DeclareMathOperator{\Tor}{Tor}

\newcommand{\co}{\colon}
\newcommand{\ra}{\rightarrow}
\newcommand{\la}{\leftarrow}
\newcommand{\lra}{\longrightarrow}

\newcommand{\thr}{\twoheadrightarrow}

\newcommand{\hra}{\hookrightarrow}

\newcommand{\ot}{\hspace{0.06em}{\otimes}\hspace{0.06em}}

\newcommand{\La}{\Leftarrow}

\newcommand{\Ra}{\Rightarrow}

\newcommand{\Llra}{\Longleftrightarrow}
\newcommand{\sbeq}{\subseteq}

\newcommand{\vare}{\varepsilon}

\newcommand{\vL}{\varLambda}

\newcommand{\limproj}{\varprojlim}

\newcommand{\wt}{\widetilde}
\newcommand{\QH}{{{}_{\vL\hspace{-0.1em}}\cat{H}}}
\newcommand{\He}{{{}_{\vL\hspace{-0.1em}}\cat{H}_{k}}}
\newcommand{\Ar}{{{}_{\vL\hspace{-0.1em}}\cat{A}_{k}}}

\newcommand{\QA}{{{}_{\vL\hspace{-0.08em}}\cat{A}}}
\newcommand{\kH}{{{}_{k}\cat{H}}}

\newcommand{\Alg}{\cat{Alg}}

\newcommand{\CM}{\cat{CM}}

\newcommand{\Df}{\hat{\cat{D}}{}^{\textnormal{fl}}}

\newcommand{\MCM}{\cat{MCM}}

\newcommand{\modf}{\cat{mod}{}^{\textnormal{fl}}}
\newcommand{\umod}{\underline{\cat{mod}}}

\newcommand{\hot}{{\tilde{\otimes}}}

\newcommand\bdot{\ensuremath{%
  \mathchoice%
   {\mskip\thinmuskip\lower0.2ex\hbox{\scalebox{1.5}{$\cdot$}}\mskip\thinmuskip}}%
   {\mskip\thinmuskip\lower0.2ex\hbox{\scalebox{1.5}{$\cdot$}}\mskip\thinmuskip}%        
   {\lower0.3ex\hbox{\scalebox{1.2}{$\cdot$}}}%  
   {\lower0.3ex\hbox{\scalebox{1.2}{$\cdot$}}}%
}

\renewcommand{\phi}{\varphi}
\renewcommand{\geq}{\geqslant}
\renewcommand{\leq}{\leqslant}

\newcommand{\fr}[1]{\mathfrak{{#1}}}

\newcommand{\cat}[1]{\mathsf{{#1}}}

\newcommand{\mc}[1]{{\mathcal{#1}}}
\newcommand{\mr}[1]{\mathrm{{#1}}}
\newcommand{\ms}[1]{\mathscr{{#1}}}
\newcommand{\tn}{\textnormal}

\newcommand{\wbar}[1]{\mkern 3.5mu\overline{\mkern-3mu{#1}\mkern-0.3mu}\mkern -0.4mu}

\newcommand{\der}[3]{{\Der}_{#1}({#2},{#3})}

\newcommand{\df}[2]{{\Def}_{\hspace{-0.08em}#2}^{#1}}

\newcommand{\hm}[4]{{\Hom}_{#2}^{#1}({#3},{#4})}

\newcommand{\uhm}[4]{{\uHom}_{#2}^{#1}({#3},{#4})}
\newcommand{\nd}[3]{{\End} _{#2}^{#1}({#3})}

\newcommand{\Q}{\hspace{-0.05em}\mathcal{O}\hspace{-0.05em}}
\newcommand{\quot}[2]{\Quot^{#1}_{#2}}

\newcommand{\shm}[4]{{\sHom} _{#2}^{#1}({#3},{#4})}

\newcommand{\tor}[4]{{\Tor}_{#2}^{#1}({#3},{#4})}
\newcommand{\ul}[1]{\underline{{#1}}}

\newcommand{\ulset}[2]{{\hspace{-0.1em}}^{#1}{\hspace{-0.2em}#2}}

\newcommand{\xra}[1]{\xrightarrow{{#1}}}

\newcommand{\xt}[4]{{\Ext}_{#2}^{#1}({#3},{#4})}

\newcommand{\uxt}[4]{{\uExt}_{#2}^{#1}({#3},{#4})}
\newcommand{\syz}[2]{{\Syz}_{#2}^{#1}}
\begin{document}
\title[Deformation theory of Cohen-Macaulay approximation]
{Deformation theory of\\ Cohen-Macaulay approximation}
\author{Runar Ile}
\address{BI Norwegian Business School, Norway}

\email{runar.ile@bi.no}
\begin{abstract}
In \cite{ile:12} we established axiomatic parametrised Cohen-Macaulay approximation which in particular was applied to pairs consisting of a finite type flat family of Cohen-Macaulay rings and modules. In this sequel we study the induced maps of deformation functors and deduce properties like smoothness and injectivity under general, mainly cohomological conditions on the module. 
\end{abstract}
\maketitle
%%%%%%%%%%%%%%%%%%%%%%%%%%%%%%%%%%%%%%%%%%%%%%%%%%%%%%%%
%%%%%%%%%%%%%%%%%%% SECTION %%%%%%%%%%%%%%%%%%%%%%%%%%%%%%%%
%%%%%%%%%%%%%%%%%%%%%%%%%%%%%%%%%%%%%%%%%%%%%%%%%%%%%%%%
\section{Introduction}\label{sec.intro}
In this article we study local properties of flat families of Cohen-Macaulay approximations by homological methods.

Let \(A\) be a Cohen-Macaulay ring of finite Krull dimension with a canonical module \(\omega_{A}\). Let \(\MCM_{A}\) and \(\cat{FID}_{A}\) denote the categories of maximal Cohen-Macaulay modules and of finite modules with finite injective dimension, respectively. M.\ Auslander and R.-O.\ Buchweitz proved in \cite{aus/buc:89} that for any finite \(A\)-module \(N\) there exists short exact sequences
\begin{equation}\label{eq.MCMseq}
0\ra L\lra M\lra N\ra 0\qquad\text{and}\qquad 0\ra N\lra L'\lra M'\ra 0
\end{equation}
with \(M\) and \(M'\) in \(\MCM_{A}\) and \(L\) and \(L'\) in \(\cat{FID}_{A}\). The maps \(M\ra N\) and \(N\ra L'\) in \eqref{eq.MCMseq} are called a maximal Cohen-Macaulay approximation and a hull of finite injective dimension, respectively, of the module \(N\).

In \cite{ile:12} we noted some of the developments since \cite{aus/buc:89}, as the study of new invariants, e.g.\
\cite{din:92, aus/din/sol:93, has/shi:97, sim/str:02} and various characterisations and applications \cite{kat:99, yos/iso:00, kat:07, her/mar:93, mar:00b, dut:04}. In his book \cite{has:00}  M.\ Hashimoto gave several new examples of the axiomatic Cohen-Macaulay approximation in \cite{aus/buc:89}. However, the `relative' and continuous aspects have received surprisingly little attention. It seems only \cite[IV 1.4.12]{has:00} and \cite{yos:99} touch upon this.

In \cite[5.1]{ile:12} we proved the following result. Let \(h\co S\ra \mc{A}\) be a Cohen-Macaulay map and \(\mc{N}\) an \(S\)-flat finite \(\mc{A}\)-module. Then there exists short exact sequences of \(S\)-flat finite \(\mc{A}\)-modules
\begin{equation}\label{eq.CMseq}
0\ra \mc{L}\lra \mc{M}\lra \mc{N}\ra 0\qquad\text{and}\qquad 0\ra \mc{N}\lra \mc{L}'\lra \mc{M}'\ra 0
\end{equation}
such that the fibres give sequences as in \eqref{eq.MCMseq} and any base change gives short exact sequences of the same kind. In the local case there are \emph{minimal} sequences \eqref{eq.CMseq} which are unique up to non-canonical isomorphisms \cite[6.2]{ile:12}. 
There are induced maps of deformation functors of pairs (algebra, module) 
\begin{equation}
\sigma_X\co \df{}{(A,N)}\lra\df{}{(A,X)}\quad\text{for}\quad X=M, M', L\,\,\text{and}\,\, L'.
\end{equation}
There are also corresponding maps \(\df{A}{N}\ra\df{A}{X}\) of the more classical deformation functors of the modules, where the algebra \(A\) only deforms trivially. To our knowledge these maps have not been defined before. They are the principal objects of study in this article. 

The main results are:
%%%%%%%%%%%%%%%%%%%%%%%%%%%
\begin{thmA*}
If\, $\xt{1}{A}{N}{M'}=0$ then $\sigma_{L'}\co \df{}{(A,N)}\lra\df{}{(A,L')}$ is formally smooth\textup{.}
If in addition $\df{}{(A,N)}$ has a versal element then so has $\df{}{(A,L')}$ and $\sigma_{L'}$ is smooth\textup{.}
\end{thmA*}
%%%%%%%%%%%%%%%%%%%%%%%%%%%
\begin{thmB*}
If\, $\xt{1}{A}{L}{N}=0$ then $\sigma_{M}\co \df{}{(A,N)}\lra\df{}{(A,M)}$ is formally smooth\textup{.}
If in addition $\df{}{(A,N)}$ has a versal element then so has $\df{}{(A,M)}$ and $\sigma_{M}$ is smooth\textup{.}
\end{thmB*}
See Theorems \ref{thm.defgrade} and \ref{thm.defgrade2} for more comprehensive statements.
The proofs of the second halves of the results use Artin's Approximation Theorem. There are analogous results for $\df{A}{N}$; see Corollaries \ref{cor.defgrade} and \ref{cor.defgrade2}.

The following is a rather surprising consequence.
%%%%%%%%%%%%%%%%%%%%%%%%%%% 
\begin{corA*}[cf. \ref{cor.defapprox}]
Each Cohen-Macaulay algebraic \(k\)-algebra \(A\) with \(A/\fr{m}_{A}\cong k\) and \(\dim A\geq 2\) has a finite \(A\)-module \(Q\) of finite \emph{projective} dimension with a \emph{universal} deformation in \(\df{A}{Q}(A)\).
\end{corA*}
Different sets of restrictions imply the main condition $\xt{1}{A}{L}{N}=0$ in Theorem B, like in the following corollary.
%%%%%%%%%%%%%%%%%%%%%%%%%%%
\begin{corB*}[cf. \ref{cor.depth}]
Assume there is a closed subscheme \(Z\) in \(\Spec A\) containing the singular locus and with complement \(U\) such that \(\tilde{N}_{\vert U}=0\) and \(\depth_{Z}N\geq 2\). Then \(\sigma_{M}\co \df{}{(A,N)}\ra\df{}{(A,M)}\) is formally smooth.
\end{corB*}
In Proposition \ref{prop.Gor} we note that $\sigma_{M}\co \df{}{(A,N)}\lra\df{}{(A,M)}$ is smooth if $A$ is Gorenstein and $\depth N=\dim A-1$, extending A.\ Ishii's \cite[3.2]{ish:00} to deformations of the pair.

Consider a quotient ring \(B=A/I\) defined by a regular sequence \(I=(f_{1},\dots,f_{n})\) and an MCM \(B\)-module \(N\). Then \(N\) is also an \(A\)-module with an MCM-approximation \(M\ra N\). Our third main result is the following (cf. Theorem \ref{thm.defMCM}):
\begin{thmC*}
Suppose $R^{N}$ and $R^M$ are minimal versal \textup{(}or formally versal\textup{)} base rings for $\df{B}{N}$ and $\df{A}{M}$. If \(N\) has a lifting to \(A/I^{2}\)\textup{,} then $R^{N}\cong R^{M}\hspace{-0.2em}/J$ for an ideal $J$ generated by linear forms.
\end{thmC*}
The proof of Theorem C is not invoking Theorems \ref{thm.defgrade} or \ref{thm.defgrade2} and applies two results which might have some independent interest. A general result applying a lifting argument gives the ideal $J$ generated by linear forms; see Lemma \ref{lem.res} (which seems to need a separability condition). 
Proposition \ref{prop.obssplit} says that the lifting condition is equivalent to the splitting of \(B\ot_{A}M\ra N\). This generalises \cite[4.5]{aus/din/sol:93} by Auslander, S. Ding, and {\O}. Solberg. The final argument shows how this splitting implies the essential conditions in Lemma \ref{lem.res}.

Theorem \ref{thm.defMCM} is illustrated by an application to hypersurface singularities where MCM-approximation is given by a functor of H. Knörrer; see Corollary \ref{cor.defMCM}.

The broader context of these results is the study of singularities in terms of the representation theory. In recent years M. Wemyss and collaborators have given many interesting results concerning the geometric McKay-Wunram correspondence where (certain) indecomposable MCM-modules correspond to irreducible components of the minimal resolution of a rational surface singularity; e.g. \cite{iya/wem:10,iya/wem:11, wem:11}. M. Van den Bergh's use in \cite{vdber:04} of endomorphism rings to prove derived equivalences for flops created a lot of activity; see Wemyss \cite{wem:18} for results and references. A general hypersurface section of certain $3$-dimensional flops gives $1$-parameter deformations of pairs (RDP-singularity, partial resolution). In \cite{gus/ile:18} we study deformation theory of pairs (rational surface singularity, MCM-module) and show that the flops are obtained by blowing up the parametrised singularity in a parametrised module. This is applied to prove several conjectures of C. Curto and D. Morrison \cite{cur/mor:13} regarding flops. The results in this article and the companion \cite{ile:19bX} contribute to the versatility of the deformation theory of pairs (algebra, module). 

The content is ordered as follows. In Section \ref{sec.CMapprox} we define the cofibred categories and in Section \ref{sec.def} the maps $\sigma_X$. We give some relevant obstruction theory for deforming modules in Section \ref{sec.obs}. The main results about the maps $\sigma_X$ with some general consequences are found in Section \ref{sec.sigma}. Section \ref{sec.defCM} concludes after several auxiliary technical results with the proof of Theorem \ref{thm.defMCM}.

Many results have analogous parts with similar arguments and the policy has been to give a fairly detailed proof of one case and leave the other cases to the reader. All rings are commutative with $1$-element. Subcategories are usually full and essential.

\subsection*{Acknowledgement}
The author thanks the referee for a detailed and helpful report.
%%%%%%%%%%%%%%%%%%%%%%%%%%%%%%%%%%%%%%%%%%%%%%%%%%%%%%%%
%%%%%%%%%%%%%%%%%%% SECTION %%%%%%%%%%%%%%%%%%%%%%%%%%%%%%%%
%%%%%%%%%%%%%%%%%%%%%%%%%%%%%%%%%%%%%%%%%%%%%%%%%%%%%%%%
\section{Preliminaries}\label{sec.CMapprox}
%%%%%%%%%%%% SUBSECTION %%%%%%%%%%%%%
\subsection{The base category}\label{subsec.base}
Fix a finite ring map \(\vL\ra k\) where $\vL$ is assumed to be excellent (in particular noetherian) \cite[\href{https://stacks.math.columbia.edu/tag/07QS}{Tag 07QS}]{SP} and $k$ a field. The kernel of $\vL\ra k$ is denoted $\fr{m}_{\vL}$. Put $k_0=\vL/\fr{m}_{\vL}$. Define $\He$ to be the category of surjective maps of $\vL$-algebras $S\ra k$ where $S$ is a noetherian, henselian, local ring \cite[\href{https://stacks.math.columbia.edu/tag/04GE}{Tag 04GE}]{SP}. A morphism is a local ring map of $\vL$-algebras $S_1\ra S_2$  commuting with the given maps to $k$. 
%%%%%%%%%%%% SUBSECTION %%%%%%%%%%%%%
\subsection{Cofibred categories}\label{subsec.loc} 
A map \(h\co S\ra \mc{A}\) of local henselian rings is \emph{algebraic} if $h$ factors as $S\ra \mc{A}^{\tn{ft}}\ra \mc{A}$ where the first map is of finite type and the second is the henselisation in a maximal ideal. Define $\Alg$ to be the category where an object is an object $S\ra k$ in $\He$ together with a map of local henselian rings $S\ra \mc{A}$ which is flat and algebraic. A morphism $(S_1\ra \mc{A}_1)\ra(S_2\ra \mc{A}_2)$ is a morphism $g\co S_1\ra S_2$ in $\He$ together with a local $S_1$-algebra map $f\co \mc{A}_1\ra \mc{A}_2$ such that the resulting commutative square is cocartesian. The fibre sum is given by the henselisation of the tensor product \(\mc{A}=\mc{A}_{1}\ot_{S_{1}}S_{2}\) in the maximal ideal \(\fr{m}_{\mc{A}_{1}}\mc{A}+\fr{m}_{S_{2}}\mc{A}\), denoted by \(\mc{A}_{1}\hot_{S_{1}}S_{2}\) or by $(\mc{A}_{1})_{S_{2}}$. It has the same closed fibre as \(S_{1}\ra \mc{A}_{1}\) and it follows that the forgetful $\Alg\ra \He$ is a cofibred category\footnote{A fibred category mimics pull-backs. We  work with rings instead of (affine) schemes. A cofibred category is a functor \(p:\cat{F}\ra\cat{C}\) such that the functor of opposite categories \(p^{\text{op}}:\cat{F}^{\text{op}}\ra\cat{C}^{\text{op}}\) is a fibred category as defined in A.\ Vistoli's \cite{vis:05}.} cofibred in groupoids; cf. \cite[\href{https://stacks.math.columbia.edu/tag/06GA}{Tag 06GA}]{SP}. In general a flat and local ring map $S\ra \mc{A}$ will be called Cohen-Macaulay if $\mc{A}\ot_SS/\fr{m}_S$ is a Cohen-Macaulay ring. There is a subcategory $\CM$ of objects in $\Alg$ which are Cohen-Macaulay maps. The forgetful $\CM\ra \He$ is a cofibred category. 

Let \(\cat{mod}\) denote the category of pairs \((h\co S\ra \mc{A},\mc{N})\) with \(h\) in \(\Alg\) and \(\mc{N}\) a finite \(\mc{A}\)-module. A morphism \((h_{1}\co S_{1}\ra \mc{A}_{1},\mc{N}_{1})\ra (h_{2}\co S_{2}\ra \mc{A}_{2},\mc{N}_{2})\) in \(\cat{mod}\) is a morphism \((g\co S_1\ra S_2,f\co \mc{A}_1\ra \mc{A}_2)\) in \(\Alg\) and an \(f\)-linear map of modules \(\alpha\co \mc{N}_{1}\ra\mc{N}_{2}\). Let $(\mc{N}_1)_{S_2}$ denote the base change $\mc{A}_2\ot_{\mc{A}_1}\mc{N}_1$. The forgetful functor \(\cat{mod}\ra\He\) is a cofibred category.
There is a cofibred subcategory \(\modf\sbeq\cat{mod}\) of modules flat over the base and a cofibred subcategory \(\MCM\sbeq\modf\) where $S\ra \mc{A}$ is in $\CM$ and \(\mc{N}\ot_{S}k\) is a maximal Cohen-Macaulay \(\mc{A}\ot_S k\)-module.

Any object \(h\co S\ra \mc{A}\) in $\CM$ has a dualising module \(\omega_{h}\) obtained by base change from the dualising module \(\omega_{h^{\tn{ft}}}\) as defined in \cite[Sec.\ 3.5]{con:00} where \(h^{\tn{ft}}\) is a finite type representative for \(h\). In particular, $(h,\omega_h)$ is an object in $\MCM$. Two finite type representatives for \(h\) factor through a common \'etale neighbourhood which is Cohen-Macaulay relative to $S$. The dualising module commutes with base change for finite type CM maps and so does \(\omega_{h}\). Let \(\cat{D}\) denote the  subcategory of \(\MCM\) of objects \((h,D)\) with \(D\) in \(\Add\{\omega_{h}\}\) and \(\Df\) the  subcategory of \(\modf\) of objects \((h\co S\ra \mc{A},\mc{N})\) such that $h$ is in $\CM$ and \(\mc{N}\) has a finite resolution by modules in \(\Add\{\omega_{h}\}\). The forgetful maps make \(\cat{D}\) and \(\Df\) cofibred categories over $\He$.

There is also a version for a fixed flat algebra. With $\vL\ra k$ as above, fix a flat ring map $\vL\ra \mc{A}$ which is the composition of two ring maps $\vL\ra \mc{A}^{\tn{ft}}\ra \mc{A}$ where the first is of finite type and the second is the henselisation at a maximal ideal. We call such an $\mc{A}$ a flat and algebraic $\vL$-algebra. There is a section $\He\ra \Alg$ given by  
$S\mapsto (S\ra \mc{A}_{S})$ where $\mc{A}_{S}=\mc{A}\hot_{\vL} S$. 
Let $\Alg^{\mc{A}}$ denote the resulting cofibred subcategory of \(\Alg\) and \(\cat{mod}_{\mc{A}}\), respectively \(\modf_{\mc{A}}\), the restriction of the cofibred categories \(\cat{mod}\) and $\modf$ to \(\Alg^{\mc{A}}\). Put $A=\mc{A}\ot_{\hspace{-0.12em}\vL}k$. If $A$ is Cohen-Macaulay then the section $\He\ra \Alg$ factors through $\CM$. Let \(\MCM_{\mc{A}}\) denote the induced cofibred subcategory of $\MCM$.
%%%%%%%%%%%%%% SUBSECTION %%%%%%%%%%%%%%%
\subsection{Deformation functors}\label{sbsec.def}
If $\vL\ra S\ra k$ is an object in $\He$ we define $\QH_{\hspace{-0.06em}S}$ as the comma category $\He/\hspace{-0.1em}S$ of maps to $S$ in $\He$.
If $h\co S\ra \mc{A}$ is an element in $\Alg$ we define $\cat{Def}_{\mc{A}/S}$ as the comma category $\Alg/(S\ra \mc{A})$, i.e. the objects are maps $(S'\ra \mc{A}')\ra(S\ra \mc{A})$ in $\Alg$ and morphism are morphisms in $\Alg$ commuting with the maps to $S\ra \mc{A}$. The objects in $\cat{Def}_{\mc{A}/S}$ are called deformations of $\mc{A}$.
If $a=(h\co S\ra \mc{A},\mc{N})$ is an object in $\modf$, we define a deformation of $a$ as a cocartesian map $a'\ra a$ in $\modf$, i.e. commutative diagram
\begin{equation}\label{eq.modfl}
\xymatrix@C-12pt@R-8pt@H-30pt{
a'\co & S'\ar[rr]^{h'}\ar[d]_(0.4)g&& \mc{A}'\ar[d]^(0.4)f & \mc{N}'\ar[d]^(0.4){\alpha} \\   
a\co & S\ar[rr]^{h}&& \mc{A} & \mc{N}
}
\end{equation}
where $(g,f)$ is an object in $\cat{Def}_{\mc{A}/S}$ (in particular $\mc{A}'_S\cong \mc{A}$) and $\alpha$ is an $f$-linear map of finite modules which are flat over the bases $S'$ and $S$, respectively, with $\mc{N}'_S\cong \mc{N}$.
A map of deformations of $a$ is a cocartesian map in $\modf$ commuting with the maps to $a$. 
Let $\cat{Def}_{(\mc{A}/S,\hspace{0.1em}\mc{N})}$ denote the resulting category of deformations of $a$. The forgetful functors make $\cat{Def}_{\mc{A}/S}$ and $\cat{Def}_{(\mc{A}/S,\hspace{0.1em}\mc{N})}$ categories cofibred in groupoids over $\QH_{\hspace{-0.06em}S}$. There is also a forgetful map of cofibred categories $\cat{Def}_{(\mc{A}/S,\hspace{0.1em}\mc{N})}\ra \cat{Def}_{\mc{A}/S}$. Similarly, fix a flat and algebraic $\vL$-algebra $\mc{A}$ as in the previous subsection, an object $S\ra k$ in $\He$ and an $S$-flat, finite $\mc{A}_S$-module $\mc{N}$. Then $\cat{Def}^{\mc{A}}_{\mc{N}\hspace{-0.15em}/\hspace{-0.1em}S}$ denotes the restriction of $\cat{Def}_{(\mc{A}_S/S,\hspace{0.1em}\mc{N})}$ to $\Alg^{\mc{A}}\hspace{-0.2em}/(S\ra \mc{A}_S)$. A deformation of $a=(S\ra \mc{A}_S,\mc{N})$ is a diagram like \eqref{eq.modfl} with $h'$ given as $(\vL\ra\mc{A})_{S'}=(S'\ra \mc{A}_{S'})$. The only possible variation is in $\mc{N}'$.

Let the deformation functors \(\df{}{\mc{A}/S}\), \(\df{}{(\mc{A}/S,\hspace{0.1em}\mc{N})}\) and \(\df{\mc{A}}{\mc{N}\hspace{-0.15em}/\hspace{-0.1em}S}\) from $\QH_{\hspace{-0.06em}S}$ to $\Sets$ be the functors corresponding to the associated groupoids of sets obtained by identifying all isomorphic objects in the fibre categories and identifying arrows accordingly. If $S=k$ with $\mc{A}=A$ and $\mc{N}=N$ we write $\df{}{A}$ for $\df{}{A/k}$ and $\df{}{(A,N)}$ for $\df{}{(A/k,N)}$. They are functors from $\He$ to $\Sets$. But with a fixed flat and algebraic $\vL$-algebra $\mc{A}$ we write $\df{\mc{A}}{N}$ for $\df{\mc{A}}{N/k}\co \He\ra \Sets$. In the case $\vL=k$ this is the classical $\df{A}{N}$.  
%%%%%%%%%%%%%% SUBSECTION %%%%%%%%%%%%%%%
\subsection{Linear approximation}
A proof of the following known result is provided in \cite[6.1]{ile:12}.
%%%%%%%%%%%% LEMMA %%%%%%%%%%%%%%%%
\begin{lem}\label{lem.Aapprox}%2
Let \(S\ra \mc{A}\) be a homomorphism of noetherian rings and \(\fr{a}\) an ideal in \(S\) such that \(I=\fr{a}\mc{A}\) is contained in the Jacobson radical of \(\mc{A}\)\textup{.} Let \(M\) and \(N\) be finite \(\mc{A}\)-modules\textup{.} Let \(\mc{A}_{n}=\mc{A}/I^{n+1}\)\textup{,} \(M_{n}=\mc{A}_{n}\ot M\) and \(N_{n}=\mc{A}_{n}\ot N\)\textup{.} Suppose there exists a tower of surjections \(\{\phi_{n}\co M_{n}\ra N_{n}\}\)\textup{.} Fix any non-negative integer \(n_{0}\)\textup{.} Then there exists an \(\mc{A}\)-linear surjection \(\psi\co M\ra N\) such that \(\mc{A}_{n_{0}}\ot\psi=\phi_{n_{0}}\)\textup{.} If the \(\phi_{n}\) are isomorphisms and \(N\) is \(S\)-flat then \(\psi\) is an isomorphism\textup{.}
\end{lem}
%%%%%%%%%%%%%%%%%%%%%%%%%%%%%%%%%%%%%%%
%%%%%%%%%%%%%%%%%% SECTION %%%%%%%%%%%%%%%%
%%%%%%%%%%%%%%%%%%%%%%%%%%%%%%%%%%%%%%%
\section{Cohen-Macaulay approximation of deformations}\label{sec.def}
We extend the Cohen-Macaulay approximation over henselian local base rings given in \cite[5.7]{ile:12} to deformations.

For each object \(a_{v}=(h_{v}\co S_{v}\ra \mc{A}_{v},\mc{N}_{v})\) in \(\modf\) with \(h_{v}\) in \(\CM\) we fix a minimal \(\MCM\)-approximation and a minimal \(\Df\)-hull
\begin{equation}\label{eq.ny} 
\pi_{v}\co 0\ra\mc{L}_{v}\ra \mc{M}_{v}\xra{\pi_{v}}\mc{N}_{v}\ra0\text{ and }
\iota_{v}\co 0\ra\mc{N}_{v}\xra{\iota_{v}}\mc{L}'_{v}\ra \mc{M}_{v}'\ra0
\end{equation}
which exist by \cite[5.7, 6.3]{ile:12}. We fix one $a=(h\co S\ra \mc{A},\mc{N})$ in \(\modf\) with \(h\) in \(\CM\), with minimal \(\MCM\)-approximation $\pi$ and minimal \(\Df\)-hull $\iota$. For each deformation \(a_{v}\ra a\), see \eqref{eq.modfl}, we choose extensions to commutative diagrams of deformations (which are all over the same deformation of algebras $h_v\ra h$) 
\begin{equation}\label{eq.right2}
\xymatrix@C-0pt@R-8pt@H-30pt{
\mc{L}_{v}\ar[r]\ar@{.>}[d]_(0.4){\lambda} & \mc{M}_{v}\ar[r]^{\pi_{v}}\ar@{.>}[d]_(0.4){\mu} & \mc{N}_{v}\ar[d]_(0.4){\nu} \\   
\mc{L}\ar[r] & \mc{M} \ar[r]^{\pi} & \mc{N}
}
\qquad\text{and}\qquad
\xymatrix@C-0pt@R-8pt@H-30pt{
\mc{N}_{v}\ar[r]^{\iota_{v}}\ar[d]_(0.4){\nu} & \mc{L}_{v}'\ar[r]\ar@{.>}[d]_(0.4){\lambda'} & \mc{M}_{v}'\ar@{.>}[d]_(0.4){\mu'} \\   
\mc{N}\ar[r]^{\iota} & \mc{L}' \ar[r] & \mc{M}'
}
\end{equation}
as follows: By \cite[6.3]{ile:12} a base change of \(\pi_{v}\) by \(S_{v}\ra S\) gives a minimal \(\MCM\)-approximation \((\mc{M}_{v})_{S}\ra(\mc{N}_{v})_S\cong \mc{N}\). By minimality it is isomorphic to \(\pi\). Choose an isomorphism. Let \(\mu\) be the composition \(\mc{M}_{v}\ra (\mc{M}_{v})_{S}\cong\mc{M}\). It is cocartesian. Do similarly for the \(\Df\)-hull. Let these choices be fixed.
%%%%%%%%%%%%%% DEFINITION %%%%%%%%%%%%%%%%%%%%
\begin{lem}\label{lem.defmap}
There are four maps
\begin{equation*}
\sigma_{X}\co \df{}{(\mc{A}/S,\hspace{0.06em}\mc{N})}\lra\df{}{(\mc{A}/S,X)}\text{ of functors }\QH_{\hspace{-0.06em}S}\lra\Sets
\end{equation*}
where \(X\) can be \(\mc{M},\mc{L},\mc{L}'\) and \(\mc{M}'\) given by \([(h_{v}\ra h, \nu)]\mapsto [(h_{v}\ra h, x)]\) for \(x\) equal to \(\mu,\lambda,\lambda'\) and \(\mu'\) in \eqref{eq.right2} respectively\textup{.}

For a flat and algebraic $\vL$-algebra $\mc{A}$ the same formulas induce well-defined maps of deformation functors of $\mc{A}$-modules
\begin{equation*}
\sigma^{\mc{A}}_{X}\co\df{\mc{A}}{\mc{N}}\lra \df{\mc{A}}{X}\,.
\end{equation*}
\end{lem}
The following lemma implies that these maps are well defined and independent of choices and thus proves Lemma \ref{lem.defmap}.
%%%%%%%%%%%%%% LEMMA %%%%%%%%%%%%%%%%%%%%%%
\begin{lem}\label{lem.defCMapprox2}
Given two deformations 
\begin{equation*}
((g_{j},f_{j}),\nu_{j})\co (h_{v_{j}}\co S_{v_{j}}\ra \mc{A}_{v_{j}},\mc{N}_{v_{j}})\ra(h_{{j}}\co S_{{j}}\ra \mc{A}_{{j}},\mc{N}_{{j}}),\quad j=1,2,
\end{equation*}
in \(\modf\) over \(\CM\)\textup{.} Consider the minimal \(\MCM\)-approximations \(\pi_{v_{j}}\) and \(\pi_{j}\) \textup{(}respectively the \(\Df\)-hulls \(\iota_{v_{j}}\) and \(\iota_{j}\)\textup{)} defined in \eqref{eq.ny} and the corresponding maps of short exact sequences \(\pi_{v_{j}}\ra\pi_{j}\) \textup{(}respectively \(\iota_{v_{j}}\ra \iota_{j}\)\textup{)} which extends \(\nu_{j}\) defined in \eqref{eq.right2}\textup{.} Given 
\begin{itemize}
\item a map \((g,f)\co h_{1}\ra h_{2}\) in $\CM$ and an \(f\)-linear map \(\alpha\co \mc{N}_{1}\ra\mc{N}_{2}\)\textup{,} 
\item maps of short exact sequences \(\pi_{1}\ra\pi_{2}\) and \(\iota_{1}\ra\iota_{2}\) which extends \(\alpha\)\textup{,}
\item a map \((\tilde{g},\tilde{f})\co h_{v_{1}}\ra h_{v_{2}}\) in $\CM$ which lifts \((g,f)\)\textup{,} and
\item an \(\tilde{f}\)-linear map \(\tilde{\alpha}\co \mc{N}_{v_{1}}\ra\mc{N}_{v_{2}}\) which lifts \(\alpha\)\textup{.}
\end{itemize}
In particular the following two diagrams of solid arrows are commutative\textup{:}
\begin{equation*}\label{eq.C}
\xymatrix@C-20pt@R-8pt@H-30pt{
& \mc{L}_{v_{2}}\ar[rr]\ar[dd]^(0.3){\lambda_{2}}|!{[dl];[dr]}\hole && \mc{M}_{v_{2}}\ar[rr]^(0.40){\pi_{v_{2}}}\ar[dd]^(0.3){\mu_{2}}|!{[dl];[dr]}\hole
&& \mc{N}_{v_{2}}\ar[dd]^(0.3){\nu_{2}} \\ 
\mc{L}_{v_{1}}\ar@{.>}[ur]\ar[rr]^(0.2){}\ar[dd]^(0.3){\lambda_{1}} && \mc{M}_{v_{1}}\ar@{.>}[ur]^(0.4){\gamma}\ar[rr]^(0.4){\pi_{v_{1}}}\ar[dd]^(0.3){\mu_{1}} 
&& \mc{N}_{v_{1}}\ar[ur]^(0.45){\tilde\alpha}\ar[dd]^(0.3){\nu_{1}} &  \\  
& \mc{L}_{2}\ar[rr]|!{[ur];[dr]}\hole && \mc{M}_{2} \ar[rr]^(0.35){\pi_{2}}|!{[ur];[dr]}\hole && \mc{N}_{2} \\
\mc{L}_{1}\ar[ur]\ar[rr] && \mc{M}_{1}\ar[ur]^(0.45){\beta\!\!} \ar[rr]^(0.45){\pi_{1}} && \mc{N}_{1}\ar[ur]_(0.47){\alpha} & 
}
\quad
\xymatrix@C-20pt@R-10pt@H-30pt{
& \mc{N}_{v_{2}}\ar[rr]^(0.40){\iota_{v_{2}}}\ar[dd]^(0.3){\nu_{2}}|!{[dl];[dr]}\hole && \mc{L}'_{v_{2}}\ar[rr]\ar[dd]^(0.3){\lambda'_{2}}|!{[dl];[dr]}\hole 
&& \mc{M}'_{v_{2}}\ar[dd]^(0.3){\mu'_{2}} \\ 
\mc{N}_{v_{1}}\ar[ur]^(0.45){\tilde\alpha}\ar[rr]^(0.35){\iota_{v_{1}}}\ar[dd]^(0.3){\nu_{1}} && \mc{L}'_{v_{1}}\ar@{.>}[ur]^(0.4){\gamma'}\ar[rr]^(0.2){}\ar[dd]^(0.3){\lambda'_{1}} 
&& \mc{M}'_{v_{1}}\ar@{.>}[ur]\ar[dd]^(0.3){\mu'_{1}} &  \\  
& \mc{N}_{2}\ar[rr]^(0.35){\iota_{2}}|!{[ur];[dr]}\hole && \mc{L}'_{2} \ar[rr]^(0.2){}|!{[ur];[dr]}\hole && \mc{M}'_{2} \\
\mc{N}_{1}\ar[ur]^(0.53){\alpha}\ar[rr]^(0.45){\iota_{1}} && \mc{L}'_{1}\ar[ur]_(0.59){\!\beta'} \ar[rr] && \mc{M}'_{1}\ar[ur] & 
}
\end{equation*}
Then there exists \(\tilde{f}\)-linear maps \(\gamma\co \mc{M}_{v_{1}}\ra\mc{M}_{v_{2}}\) and \(\gamma'\co \mc{L}'_{v_{1}}\ra\mc{L}'_{v_{2}}\) such that the induced diagrams are commutative\textup{.} If \(\tilde\alpha\) is cocartesian\textup{,} so are \(\gamma\) and \(\gamma'\)\textup{.}
\end{lem}
\begin{proof}
Consider the \(\MCM\)-approximation case. By applying base changes to the front diagram, we can reduce the problem to the case \(h_{v_{1}}\ra h_{1}\) equals \(h_{v_{2}}\ra h_{2}\). Then, by \cite[5.7]{ile:12}, there is a lifting \(\gamma_{1}\co \mc{M}_{v_{1}}\ra\mc{M}_{v_{2}}\) of \(\tilde\alpha\). We would like to adjust \(\gamma_{1}\) so that it lifts \(\beta\) too. We have that \(\mu_{2}\gamma_{1}-\beta\mu_{1}\) factors through \(\mc{L}_{2}\) by a map \(\tau\co \mc{M}_{v_{1}}\ra\mc{L}_{2}\). It induces a unique map \(\wbar\tau\co \mc{M}_{1}\ra\mc{L}_{2}\) since \(\mu_{1}\) is cocartesian. If \(\mc{D}_{*}\thr\mc{L}_{v_{2}}\) is a finite \(\cat{D}\)-resolution, then base change gives a finite \(\cat{D}\)-resolution \(\mc{D}_{*}\ot_{S_{v_{2}}}S_{2}\thr \mc{L}_{2}\) and \(\wbar\tau\) factors through a \(\wbar\sigma\co \mc{M}_{1}\ra\mc{D}_{0}\ot_{S_{v_{2}}}S_{2}\) by \cite[5.7]{ile:12}. Since $\hm{}{}{\mc{M}_{v_{1}}}{\mc{D}_{0}}$ is a deformation of $\hm{}{}{\mc{M}_{1}}{\mc{D}_{0}\ot_{S_{v_{2}}}S_{2}}$ (cf. \cite[2.4]{ile:12})
there is a \(\sigma\co \mc{M}_{v_{1}}\ra\mc{D}_{0}\) lifting \(\wbar\sigma\). Subtracting the induced map \(\mc{M}_{v_{1}}\ra\mc{M}_{v_{2}}\) from \(\gamma_{1}\) gives our desired \(\gamma\). If \(\tilde\alpha\) is an isomorphism so is \(\gamma\) by minimality of the approximations \(\pi_{v_{j}}\).
The argument for the \(\Df\)-case is similar.
\end{proof}
%%%%%%%%%%%%%%%%%%%%%%%%%%%%%%%%%%%%%%%%%%%%%%
%%%%%%%%%%%%%%%%%%%% SECTION %%%%%%%%%%%%%%%%%%%%%
%%%%%%%%%%%%%%%%%%%%%%%%%%%%%%%%%%%%%%%%%%%%%%
\section{Obstructions}\label{sec.obs}
We summarise some obstruction theory for deformations of modules which will be used to study the maps in Lemma \ref{lem.defmap}.

Suppose $\beta\co 0\ra J\xra{\;i\,} B'\xra{\;q\;} B\ra 0$ is an extension of rings where $i$ denotes the inclusion of the ideal $J$. Assume $J^2=0$. If $N'$ is a $B'$-module, $N'\ot_{B'}-$ applied to $\beta$ gives the exact sequence $0\ra\tor{B'}{1}{N'}{B}\ra N'\ot_{B'}J\xra{{\id}\ot i}N'\ra N'\ot_{B'}B\ra 0$. Note that $J$ is a $B$-module as $B'$-module since $J^2=0$. It follows that $N'\ot_{B'}J\cong (N'\ot_{B'}B)\ot_BJ$. 
%%%%%%%%% DEFINITION %%%%%%%%%%%%%
\begin{defn}\label{def.modulobs}
Given an extension $\beta$ with $J^2=0$ as above and suppose $N$ is a $B$-module. Then a $B'$-module $N'$ with a surjection of $B'$-modules 
$\alpha\co N'\ra N$ is called a \emph{lifting} of $N$ along $q$ (or to $B'$) if the natural $B'$-linear map $j\co N\ot_BJ\ra N'$ defined by $j(n\ot u)=u\cdot\tilde{n}$ for any $\tilde{n}\in N'$ with $\alpha(\tilde{n})=n$ gives an isomorphism $N\ot_BJ\cong \ker(\alpha)$. Two liftings $\alpha_i\co N'_i\ra N$ ($i=1,2$) along $q$ are \emph{equivalent} if there is an isomorphism $\phi\co N_1'\ra N_2'$ of $B'$-modules with $\alpha_1=\alpha_2\phi$.
\end{defn}
For a lifting $\alpha$ it follows that $N'\ot_{B'} B\cong N$ and $\tor{B'}{1}{N'}{B}=0$, and vice versa, a surjection $\alpha$ is a lifting if these two conditions hold. Moreover, $\alpha$ gives a $B'$-module extension $\nu\co 0\ra N\ot J\xra{\;j\;} N'\xra{\;\alpha\;} N\ra 0$. Two liftings along $q$ are equivalent if and only if the corresponding extensions are isomorphic.
There is an obstruction theory for liftings of modules in terms of \(\Ext\) groups.
%%%%%%%%%%%%%%% PROPOSITION %%%%%%%%%%%%%%%%%
\begin{prop}\label{prop.obsmodule}
Given an extension $\beta$ and a $B$-module $N$ as in \textup{Definition \ref{def.modulobs}}\textup{.}
\begin{enumerate}[leftmargin=2.4em, label=\textup{(\roman*)}]
\item There is an element \(\ob(q,N)\in\xt{2}{B}{N}{N\ot_B J}\)  
such that \(\ob(q,N)=0\) if and only if there exists a lifting of \(N\) along \(q\)\textup{.}
\item If \(\ob(q,N)=0\) then the set of equivalence classes of liftings of \(N\) along \(q\) is a torsor for \(\xt{1}{B}{N}{N\ot_B J}\)\textup{.} 
\item The set of automorphisms of a given lifting is canonically isomorphic to $\hm{}{B}{N}{N\ot_B J}$\textup{.}
\end{enumerate}
\end{prop}
\begin{proof}
(i) Pick a $B$-free resolution of $N$: $\cdots F_2\xra{d_2}F_1\xra{d_1}F_0\xra{\vare}N\ra 0$. Lift the differentials to maps $\tilde{d}_n\co \tilde{F}_n\ra \tilde{F}_{n-1}$ of $B'$-free modules of the same rank. Denote the map $\tilde{F}\ra F$ by $\pi$. Then $\tilde{d}_1\tilde{d}_2$ is induced by a map $\eta_2\co F_2\ra F_0\ot J$ and $\wbar{\eta}_2:=(\vare\ot \id)\eta_{2}$ is a $2$-cocycle in the complex $\hm{}{B}{F}{N\ot J}$ since $\eta_2d_3$ equals $(d_1\ot\id_J)\eta_3$ where $\eta_3\co F_3\ra F_1\ot J $ is inducing $\tilde{d}_2\tilde{d}_3$. The class of $\wbar{\eta}_2$ defines $\ob(q,N)$. It is independent of the chosen resolution and liftings. 

If there is a lifting $N'$ of $N$ along $q$ we can choose a $B'$-free resolution $F'$ of $N'$ first. Then $\cH^0(B\ot F')\cong N$ and $\cH^1(B\ot F')=0$. A $B$-free resolution $F$ of $N$ is obtained by adding terms in degree $\geq 3$. It follows that $\ob(q,N)=0$. 

Suppose $\ob(q,N)=0$. Then there is a $\wbar{\xi}\co F_1\ra N\ot J$ with $\wbar{\eta}_2=\wbar{\xi} \pi$. Let $\xi_1\co F_1\ra F_0\ot J$ be a lifting of $\wbar{\xi}$ and let $\iota$ denote the inclusion $\iota\co F\ot J\ra \tilde{F}$. Let $\xi_2\co F_2\ra F_1\ot J$ be a lifting of $\eta_{2}-\xi_1d_2$. Then $(\tilde{d}_1-\iota\xi_1 \pi)(\tilde{d}_2-\iota\xi_2 \pi)=0$ which implies that $N':=\coker(\tilde{d}_1-\iota\xi_1 \pi)$ with its natural map to $N$ gives a lifting of $N$ along $q$.

(ii) Given two liftings $N'_1$ and $N'_2$ of $N$ along $q$. By what we did above there are maps $d'_{i,1}, d'_{i,2}\co F'_i\ra F'_{i-1}$ for $i=1,2$ such that $B\ot_{B'}d'_{i,j}=d_i$. Then $d_{i,2}-d_{i,1}$ equals $\iota\xi_i\pi$ for some $\xi_i$. One calculates
\begin{equation}
\iota(\xi_1d_2+(d_1\ot\id)\xi_2)\pi=(\iota\xi_1\pi)d_{2,2}+d_{1,1}(\iota\xi_2\pi)=0
\end{equation}
which implies that $\xi_1d_2+(d_1\ot\id)\xi_2=0$. Then $\wbar{\xi}:=(\vare\ot\id)\xi_1$ defines a class in $\xt{1}{B}{N}{N\ot J}$. Conversely, given a lifting $N'$ with differential $d'$, such a class can be lifted to maps $\xi_1$ and $\xi_2$ with $\xi_1d_2+(d_1\ot\id)\xi_2=0$. Then $\coker(d'_1+\iota\xi_1\pi)$ defines another lifting.

(iii) follows since automorphisms of the lifting $\alpha$ equals automorphisms of the corresponding extension $\nu$ and $\id_{N'}$ corresponds to $0$ in $\hm{}{B}{N}{N\ot J}$.
\end{proof}
The element \(\ob(q,N)\) is called the \emph{obstruction} of \((q,N)\). If $N'\ra N$ is a lifting and $\xi$ an element in $\xt{1}{B}{N}{N\ot J}$ we write $(N'+\xi)\ra N$ for the new lifting obtained by the torsor action.
%%%%%%%%%%%%%% LEMMA %%%%%%%%%%%%%%
\begin{lem}\label{lem.obsdef}
Given a commutative diagram
\begin{equation*}
\xymatrix@C+0pt@R-12pt@H-30pt{
R'\ar[r]\ar[d]_(0.4)r & S'\ar[d]^(0.4)s
\\
R\ar[r] & S
}
\end{equation*}
in $\He$ where $r$ and $s$ are surjective. Put $H=\ker(r)$ and $I=\ker(s)$ and assume $H^2=0=I^2$\textup{.} Let $R'\ra B'$ be an object in $\cat{Alg}$, let $R\ra B=B'_R$ denote the base change to $R$ and let $q\co B'\ra B$ denote the induced map. Suppose $N$ is an $R$\textup{-}flat $B$\textup{-}module\textup{.} Then\textup{:}
\begin{enumerate}[leftmargin=2.4em, label=\textup{(\roman*)}]
\item Base change of $q$ to $S'$ gives the extension $0\ra B_S\ot_SI\ra B'_{S'}\ra B_{S}\ra 0$.
\item There is a natural isomorphism $$\xt{n}{B_{S}}{N_S}{N_S\ot_{B_S}(B_S\ot_SI)}\cong\xt{n}{B}{N}{N_S\ot_{S}I}\tn{ for all } n.$$
\item The obstruction element \(\ob(q,N)\) maps to $\ob(q_{S'},N_S)$ along the map $$\tau^2_*\co\xt{2}{B}{N}{N\ot_RH}\lra \xt{2}{B}{N}{N_S\ot_SI}$$ induced by the natural map $\tau\co N\ot_RH\ra N_S\ot_SI$\textup{.}
\item The torsor action commutes with base change\textup{:} If $N'\ra N$ is a lifting of $N$ along $q$ and $\xi\in\xt{1}{B}{N}{N\ot_R H}\) then the base change to $S'$ of the lifting $(N'+\xi)\ra N$ is equivalent to $(N'_{S'}+\tau^1_*(\xi))\ra N_S$\textup{.} 
\end{enumerate}
\end{lem}
\begin{proof}
(i) Base change of $\ker(q)\ra B'\ra B$ to $S'$ equals $B'_{S'}\ot_{S'}(I\ra S'\ra S)$ which gives a short exact sequence. Moreover, $B'_{S'}\ot_{S'}I\cong B_{S}\ot_{S}I$. (ii) follows by a change of rings spectral sequence (cf. \cite[2.3.1]{ile:12}). 

(iii) If $(F,d)\thr N$ is a $B$-free resolution of $N$, then the base change $ (F_S,d_S)\thr N_S$ is a $B_S$-free resolution of $N_S$. Then the base change $\tilde{d}_{S'}$ of a lifting $\tilde{d}$ of the differential $d$ is a lifting of the differential $d_{S'}$. The obstruction $\ob(q_{S'},N_S)$ is induced by $(\tilde{d}_{S'})^{2}$ which equals $(\tilde{d}\,^2)_{S'}$, i.e. the base change of the map which induces $\ob(q,N)$. The map in (ii) then takes $\ob(q_{S'},N_S)$ to $\tau^2_*\ob(q,N)$. (iv) is similar.
\end{proof}
If $\fr{m}_RH=0$ (i.e. $R'\ra R$ is small) then by (ii)
\begin{equation}\label{eq.xtred}
\xt{n}{B}{N}{N\ot_RH}\cong \xt{n}{B\ot_Rk}{N\ot_Rk}{N\ot_Rk}\ot_kH\,,
\end{equation}
so in this case there are fixed $k$-vector spaces which classify obstruction and give the torsor action.

If, in the setting of Lemma \ref{lem.obsdef}, $N$ is also finite (so $(B,N)$ is a deformation of $(B\ot_Rk,N\ot_Rk)$) with $\ob(q,N)=0$, then a lifting $N'$ can be chosen to be finite by the proof of Proposition \ref{prop.obsmodule}. Moreover, a $B'$-free resolution $\vare\co F\thr N'$ is also an $R'$-flat resolution of $N'$. Then $0=\tor{B'}{1}{B}{N'}\cong \cH_1(B\ot_{B'}F)\cong \cH_1(R\ot_{R'}F)\cong \tor{R'}{1}{R}{N'}$ which implies that $\tor{R'}{1}{k}{N'}\cong\cH_1(k\ot_{R'}F)\cong\cH_1(k\ot_{R}R\ot_{R'}F)\cong \tor{R}{1}{k}{N}=0$. By the local criterion of flatness, $N'$ is $R'$-flat and so $(B',N')$ is a deformation of $(B,N)$. For brevity we will also simply say that $N'$ is a deformation of $N$.
%%%%%%%%%%%%%% LEMMA %%%%%%%%%%%%%%
\begin{lem}\label{lem.obsmodule}
Given an extension $\beta$ and a $B$-module $N$ as in \textup{Definition \ref{def.modulobs}}\textup{.} Suppose $\vare\co (F,d)\thr N$ is a $B'$\textup{-}free resolution of $N$\textup{.} Put $N'_1=\ker\vare$\textup{,} $\wbar{N}{\hspace{0.1em}}'_{\hspace{-0.15em}1}=B\ot_{B'}N'_1$ and $\wbar{F}=B\ot_{B'}F$\textup{.} Applying $B\ot_{B'}-$ to the short exact sequence $0\ra N'_1\ra F_0\ra N\ra 0$ gives a $4$-term exact sequence 
\begin{equation*}
0\ra N\ot_BJ\lra \wbar{N}{\hspace{0.1em}}'_{\hspace{-0.15em}1}\lra \wbar{F}_{\hspace{-0.1em}0}\lra N\ra 0
\end{equation*} 
which represents $\ob(q,N)$ in $\xt{2}{B}{N}{N\ot J}$\textup{.}
\end{lem}
\begin{proof}
The $4$-term exact sequence is obtained since $\tor{B'}{1}{N}{B}\cong N\ot_BJ$.
Choose a surjection $\gamma\co E\ra J$ where $E$ is $B'$-free.  
Since $NJ=0$, the composition with the multiplication map $F_0\ot E\ra F_0\ot J\ra F_0$ factors through a $B'$-linear map $\psi\co F_0\ot E\ra F_1$. Then 
\begin{equation}
\xymatrix@C+6pt@R-26pt@H-30pt{
&&& \hspace{-0.45em}\wbar{F}_{\hspace{-0.1em}0}\ot \wbar{E}
\\
&&& \oplus
\\
N & \wbar{F}_{\hspace{-0.1em}0}\ar@{->>}[l]_(0.46){\wbar{\vare}} & \wbar{F}_{\hspace{-0.1em}1}\ar[l]_(0.44){\wbar{d}_{1}} & 
\wbar{F}_{\hspace{-0.1em}2}
\ar[l]_(0.46){(\wbar{\psi}\,,\,\wbar{d}_{2})}
}
\end{equation}
gives a $B$-free $2$-presentation of $N$; cf. \cite[Lemma 3]{ile:04a}. Following the proof of Proposition \ref{prop.obsmodule}, $\wbar{\eta}_2$ can be given by $\wbar{\vare}\ot\wbar{\gamma}\co \wbar{F}_0\ot \wbar{E}\ra N\ot J$. Since the upper row in the commutative diagram 
\begin{equation}
\xymatrix@C+6pt@R-26pt@H-30pt{
&&&& \hspace{-0.45em}\wbar{F}_{\hspace{-0.1em}0}\ot \wbar{E}
\\
&&&& \oplus
\\
&& \wbar{F}_{\hspace{-0.1em}0}\ar@{=}[ddddd]
& \wbar{F}_{\hspace{-0.1em}1}\ar[ddddd]\ar[l]_(0.44){\wbar{d}_{1}} & 
\wbar{F}_{\hspace{-0.1em}2}\ar[l]\ar[ddddd]^(0.44){(\wbar{\vare}\ot\wbar{\gamma}\,,\,0)}
\ar[l]_(0.46){(\wbar{\psi}\,,\,\wbar{d}_{2})} & \cdots\ar[ddddd]\ar[l]
\\
\\
\\
\\
\\
0 & N\ar[l] & \wbar{F}_{\hspace{-0.1em}0}\ar[l]_(0.46){\wbar{\vare}} & \wbar{N}{\hspace{0.1em}}'_{\hspace{-0.15em}1}\ar[l] & 
N\ot J
\ar[l] & 0\ar[l]
}
\end{equation}
is the beginning of a $B$-free resolution of $N$, $\wbar{\vare}\ot\wbar{\gamma}$ also defines the image of the $4$-term exact sequence in $\xt{2}{B}{N}{N\ot J}$.
\end{proof}
%%%%%%%%%%%%%% LEMMA %%%%%%%%%%%%%%
\begin{lem}\label{lem.omap}
Let \(k\) be a field and \(A\) a local algebraic \(k\)\textup{-}algebra\textup{.} Given a small surjection \(p\co R\ra S\) in \(\He\) and a deformation \(R\ra \mc{A}\) of \(k\ra A\)\textup{.} Put \(q=\id\!\ot 1\co \mc{A}\ra \mc{A}\ot_{R}S=\mc{A}_{S}\) and \(I=\ker p\)\textup{.} 
Given commutative diagrams 
\begin{equation*}
\xymatrix@C-6pt@R-9pt@H-30pt{
0\ar[r] & \mc{N}\ar[r]^{\iota}\ar[d] & \mc{L}'\ar[r]^{\eta}\ar[d] & \mc{M}'\ar[r]\ar[d] & 0 \ar@{}[drr]|-{\text{\normalsize and}}
&&
0\ar[r] & \mc{L}\ar[r]^{\rho}\ar[d] & \mc{M}\ar[r]^{\pi}\ar[d] & \mc{N}\ar[r]\ar[d] & 0
\\
0\ar[r] & N\ar[r]^{\wbar{\iota}} & L'\ar[r]^{\wbar{\eta}} & M'\ar[r] & 0
&&
0\ar[r] & L\ar[r]^{\wbar{\rho}} & M\ar[r]^{\wbar{\pi}} & N\ar[r] & 0
}
\end{equation*}
with short exact horizontal sequences\textup{,} the upper of \(\mc{A}_{S}\)\textup{-}modules and the lower of \(A\)\textup{-}modules, where the vertical maps are deformations.
Then\textup{:}
\begin{enumerate}[leftmargin=2.4em, label=\textup{(\roman*)}]
\item \(\wbar{\iota}_{*}\ob(q,\mc{N})=\wbar{\iota}^{*}\ob(q,\mc{L}')\) in \(\xt{2}{A}{N}{L'}\ot_{k} I\)
\item \(\wbar{\pi}^{*}\ob(q,\mc{N})=\wbar{\pi}_{*}\ob(q,\mc{M})\) in \(\xt{2}{A}{M}{N}\ot_{k}I\)
\end{enumerate}
Furthermore\textup{,} assume we have short exact sequences of \(\mc{A}\)-modules
\begin{equation*}
\xymatrix@C-6pt@R-9pt@H-30pt{
0\ar[r] & \wt{\mc{N}}_{i}\ar[r]^(0.48){\iota_i} & \wt{\mc{L}}_{i}'\ar[r]^(0.45){\eta_i} & \wt{\mc{M}}_{i}'\ar[r] & 0
}
\text{ and }
\xymatrix@C-6pt@R-9pt@H-30pt{
0\ar[r] & \wt{\mc{L}}_{i}\ar[r]^(0.45){\rho_i} & \wt{\mc{M}}_{i}\ar[r]^(0.48){\pi_i} & \wt{\mc{N}}_{i}\ar[r] & 0
}
\end{equation*}
for \(i=1,2\) with maps to the corresponding upper sequences above which form commutative diagrams of deformations\textup{.} Let \(\delta,\zeta\) and $\xi$ denote the differences of the deformations \(\wt{\mc{N}}_{i}\ra\mc{N}\)\textup{,} the \(\wt{\mc{L}}'_{i}\ra\mc{L}'\) and the \(\wt{\mc{M}}_{i}\ra\mc{M}\), respectively \textup{(}cf\textup{.}\ Proposition \ref{prop.obsmodule}\textup{).} Then\textup{:}
\begin{enumerate}[resume*]
\item \(\wbar{\iota}_{*}\delta=\wbar{\iota}^{*}\zeta\) in \(\xt{1}{A}{N}{L'}\ot_{k} I\) and \(\wbar{\pi}^{*}\delta=\wbar{\pi}_{*}\xi\) in \(\xt{1}{A}{M}{N}\ot_{k}I\)
\end{enumerate}
\end{lem}
\begin{proof}
(i) Let \((\mc{F},d)\) be an $\mc{A}_S$-free resolution of \(\mc{N}\) and put \(F:=\mc{F}\ot_{S}k\) which is an $A$-free resolution of \(N\):
\begin{equation}
\xymatrix@C-4pt@R-9pt@H-30pt{
0 & \mc{N}\ar[l]\ar[d] & \mc{A}_{S}^{n_{0}}\ar[l]_{\vare}\ar[d] & \mc{A}_{S}^{n_{1}}\ar[l]_{d_{1}}\ar[d] & \mc{A}_{S}^{n_{2}}\ar[l]_{d_{2}}\ar[d] & \dots\ar[l]_{d_{3}}
\\
0 & N\ar[l] & A^{n_{0}}\ar[l]_{\wbar{\vare}} & A^{n_{1}}\ar[l]_{\wbar{d}_{1}} & A^{n_{2}}\ar[l]_{\wbar{d}_{2}} & \dots\ar[l]_{\wbar{d}_{3}}
}
\end{equation}
Similarly, let \((\mc{G},d')\) be a free resolution of \(\mc{M}'\) and put \(G=\mc{G}\ot_{S}k\). Then one can make \(\mc{F}\oplus\mc{G}\) a free resolution of \(\mc{L}'\) with a differential of the form \(\bigl(\begin{smallmatrix}
d' & 0 \\ e & d
\end{smallmatrix}
\bigr)
\).
To find the obstruction we lift the differentials to maps of free \(\mc{A}\)-modules: \(\tilde{d}_{1}\co \mc{A}^{n_{1}}\ra \mc{A}^{n_{0}}\) lifts \(d_{1}\) and so on. Then the obstruction for lifting \(\mc{N}\) to \(\mc{A}\) is induced by \(\tilde{d}^{\hspace{0.1em}2}\) which factors through a degree two cocycle \(a\co F\ra F\ot_{k}I\) in the Yoneda complex \(\nd{\bdot}{A}{F}\ot_{k}I\) which represents \(\ob(q,\mc{N})\). In the case of \(\mc{L}'\) the obstruction is induced by 
\(\bigl(\begin{smallmatrix}
\tilde{d}' & 0 \\ \tilde{e} & \tilde{d}
\end{smallmatrix}
\bigr){}^{2}
\)
which factors through a degree two cocycle
\(\bigl(\begin{smallmatrix}
b & 0 \\ c & a
\end{smallmatrix}
\bigr)
\) in \(\nd{\bdot}{A}{G\oplus F}\ot_{k}I\) which represents \(\ob(q,\mc{L}')\). Since \(\wbar{\iota}\) is represented by the inclusion of resolutions \(F\ra G\oplus F\) we find that \(\wbar{\iota}_{*}a=
 \bigl(\begin{smallmatrix}
0 \\ a
\end{smallmatrix}
\bigr)
=
\wbar{\iota}^{*}\bigl(\begin{smallmatrix}
b & 0 \\ c & a
\end{smallmatrix}
\bigr)\) which in cohomology gives \(\wbar{\iota}_{*}\ob(q,\mc{N})=\wbar{\iota}^{*}\ob(q,\mc{L}')\). A similar argument gives (ii).

\vspace{-0.1em}
(iii), first part: We can assume that \(\wt{\mc{L}}'_{i}\) has a resolution with differential 
\(\bigl(\begin{smallmatrix}
\tilde{d}_{i}' & 0 \\
\tilde{e}_{i} & \tilde{d}_{i}
\end{smallmatrix}
\bigr)
\) for \(i=1,2\) 
lifting the resolution of \(\mc{L}'\) given above. Then the difference of the two differentials factors through a degree one cocycle
\(\bigl(\begin{smallmatrix}
s & 0 \\ t & r
\end{smallmatrix}
\bigr)
\) in \(\nd{\bdot}{A}{G\oplus F}\ot_{k}I\) which represents \(\zeta\). Then the rest is analogous to (ii). The second part is similar. 
\end{proof}
%%%%%%%%%%%%%%%%%%%%%%%%%%%%%%%%%%%%%%%%%%%%%%
%%%%%%%%%%%%%%%%%%%% SECTION %%%%%%%%%%%%%%%%%%%%%
%%%%%%%%%%%%%%%%%%%%%%%%%%%%%%%%%%%%%%%%%%%%%%
\section{Maps of deformation functors induced by\\ Cohen-Macaulay approximation}\label{sec.sigma}
After two lemmas relating to the Schlessinger-Rim conditions in Artin's \cite{art:74} we state several results about various maps of deformation functors induced by Cohen-Macaulay approximation.

For any cofibred category $\cat{F}$ over $\He$ (or over the  subcategory $\Ar$ of Artin rings) we will in the following assume that the fibre category $\cat{F}(k)$ is equivalent to a one-object, one-morphism category. Furthermore, for all maps $f\co R\ra S$ in $\He$ and for all objects $a$ in $\cat{F}(R)$ we choose a push forward $f_*a$ in $\cat{F}(S)$. Let $F=\bar{\cat{F}}$ denote the functor associated to $\cat{F}$.
%%%%%%%%% DEFINITION %%%%%%%%%%%%%
\begin{defn}\label{def.smooth}
Assume that $\cat{F}$ and $\cat{G}$ are cofibred categories over $\He$ which are locally of finite presentation (`limit preserving' in \cite[p. 167]{art:74}).
A map \(\phi\co \cat{F}\ra \cat{G}\) is \emph{smooth} (formally smooth) if, for all surjections \(f\co S'\ra S\) in \(\He\) (respectively in \(\Ar\)), the natural map 
\begin{equation}\label{eq.vers}
(f_*,\phi(S'))\co F(S')\ra F(S)\times_{G(S)}G(S')
\end{equation}
is surjective. Put $h_R=\hm{}{\He}{R}{-}$. Let $v$ be an object in \(\cat{F}(R)\) and let $c_v\co h_R\ra F$ denote the corresponding Yoneda map. 
If \(R\) is algebraic as \(\vL\)-algebra and $c_v$ is smooth (an isomorphism) then $v$ is \emph{versal} (respectively universal). Moreover, $v$ (or a formal element $v=\{v_{n}\}$ in $\limproj \cat{F}(R/\fr{m}_{R}^{n+1})$) is \emph{formally versal} if $c_v$ restricted to $\Ar$ is formally smooth.
\end{defn}
%%%%%%%%% DEFINITION %%%%%%%%%%%%%
\begin{defn}\label{def.S2}
Suppose $\cat{F}\ra \He$ is a cofibred category satisfying the Schlessinger-Rim condition (S1') in \cite[2.2]{art:74} with associated functor $F$. Let $a$ be an object in $\cat{F}(S)$ and $I$ a finite $S$-module. Put $S{\oplus}I=\Sym_S(I)/(\Sym_S^2(I))$. Let $\cat{F}_a(S{\oplus}I)$ denote the groupoid of maps $a'\ra a$ above the projection $p\co S{\oplus}I\ra S$ and let $D^{F}_a(I)$ denote the $S$-module of isomorphism classes in $\cat{F}_a(S{\oplus}I)$; \cite[2.10]{sch:68}. Define the condition on $\cat{F}$:
\begin{itemize}
\item[(S2)] $D^{F}_a(I)$ is a finite $S$-module
\end{itemize}
for all reduced $S$ in $\He$, objects $a$ and finite $S$-modules $I$; \cite[2.5]{art:74}.
\end{defn}
If \(A\) is a local algebraic \(k\)-algebra and \(N\) a finite \(A\)-module then by standard arguments $\cat{Def}_{(A,N)}$ is locally of finite presentation and satisfies (S1'); cf. \cite[4.1]{ile:19bX}, and likewise for $\cat{Def}^{\mc{A}}_{N}$ where $\mc{A}$ is a flat and algebraic $\vL$-algebra.
%%%%%%%%%%%%%% LEMMA %%%%%%%%%%%%%%
\begin{lem}\label{lem.S2v}
Suppose $\cat{F}$ satisfies \textup{(S1')} and has a versal object $v$ in $\cat{F}(R)$. Then $\cat{F}$ satisfies \textup{(S2).} 
\end{lem}
\begin{proof}
We use the assumptions in Definition \ref{def.S2}.
By versality there is a $g$ in $h_R(S)$ with $g_*v\cong a$ in $\cat{F}(S)$. If $a'$ is a lifting of $a$ along $p$ then there is a $g'$ lifting $g$ with $g'_*v\cong a'$ by versality. I.e. the $S$-linear map $D_g^{h_R}(I)\ra D_a^{F}(I)$ is surjective. Now $D_g^{h_R}(I)\cong \hm{}{R}{\Omega_{R/\vL}}{I}$. Since $R$ is algebraic, $\Omega_{R/\vL}$ is a finite $R$-module and so is $D_a^{F}(I)$.
\end{proof}
Let \(A\) be a Cohen-Macaulay local algebraic \(k\)-algebra and \(N\) a finite \(A\)-module.
Fix a minimal \(\cat{MCM}_{A}\)-approximation \(0\ra L\ra M\xra{\,\pi\,} N\ra0\) and a minimal \(\hat{\cat{D}}_{A}\)-hull \(0\ra N\xra{\,\iota\,} L'\ra M'\ra0\).
%%%%%%%%%%%%%% LEMMA %%%%%%%%%%%%%%
\begin{lem}\label{lem.S2}
Suppose \textup{(S2)} holds for \(\df{}{(A,N)}\)\textup{.} 
\begin{enumerate}[leftmargin=2.4em, label=\textup{(\roman*)}]
\item If $\xt{1}{A}{N}{M'}=0$ then \textup{(S2)} holds for $\df{}{(A,L')}$\textup{.}
\item If $\xt{1}{A}{L}{N}=0$ then \textup{(S2)} holds for $\df{}{(A,M)}$\textup{.}
\end{enumerate}
\end{lem}
\begin{proof}
(i) We use the assumptions in Definition \ref{def.S2}. 
Let $a=(S\ra\mc{A},\mc{N})\in \df{}{(A,N)}(S)$, and consider the bottom short exact sequence to the right in \eqref{eq.right2}. Let $a_0=(S\ra\mc{A})\in \df{}{A}(S)$ be the image of $a$ by the forgetful map.
Suppose $r>0$. As $\xt{r}{A}{M'}{L'}=0$, base change theory implies that $\xt{r}{\mc{A}}{\mc{M}'}{\mc{L}'\ot I}=0$; cf. \cite[5.1]{ogu/ber:72}.  Then the natural map $(\iota_{\mc{N}})^*\co\xt{r}{\mc{A}}{\mc{L}'}{\mc{L}'\ot I}\ra \xt{r}{\mc{A}}{\mc{N}}{\mc{L}'\ot I}$ is an isomorphism. Composing the surjection $(\iota_{\mc{N}})_*\co \xt{r}{\mc{A}}{\mc{N}}{\mc{N}\ot I}\ra \xt{r}{\mc{A}}{\mc{N}}{\mc{L}'\ot I}$ with the inverse of $(\iota_{\mc{N}})^*$ gives a natural map 
\begin{equation}
\eta^r\co \xt{r}{\mc{A}}{\mc{N}}{\mc{N}\ot I}\ra \xt{r}{\mc{A}}{\mc{L}'}{\mc{L}'\ot I}\,.
\end{equation}
Base change theory and the assumption implies as above that $\xt{1}{\mc{A}}{\mc{N}}{\mc{M}'\ot I}=0$. From the long exact sequence it follows that $\eta^1$ is surjective and $\eta^2$ is injective. 

Put $D_{a_0}(I)=\df{}{\mc{A}/S}(S{\oplus}I)$ and $b=(S\ra\mc{A},\hspace{0.1em}\mc{L}')$. To the ring maps $S\ra\mc{A}\ra\mc{A}{\oplus}\mc{N}$ there is a natural Jacobi-Zariski long-exact sequence of (graded) André-Quillen cohomology obtained from \cite[Chap. IV, 2.3]{ill:71} which maps to the corresponding sequence for $k\ra A\ra A{\oplus} N$, see \cite[2.10]{ile:19bX}. Low-degree terms give the commutative diagram of $\mc{A}$-modules 
\begin{equation}
\xymatrix@C-6pt@R-6pt@H-10pt{
\xt{1}{\mc{A}}{\mc{N}}{\mc{N}\ot I}\ar[r]\ar[d]^(0.43){\eta^1} & D_a(I)\ar[r]\ar[d]^(0.43){\delta} & D_{a_0}(I)\ar[r]\ar@{=}[d] & \xt{2}{\mc{A}}{\mc{N}}{\mc{N}\ot I}\ar[d]^(0.43){\eta^2}
\\
\xt{1}{\mc{A}}{\mc{L}'}{\mc{L}'\ot I}\ar[r] & D_b(I)\ar[r] & D_{a_0}(I)\ar[r] & \xt{2}{\mc{A}}{\mc{L}'}{\mc{L}'\ot I}
}
\end{equation}
where the middle terms are canonically isomorphic to the degree one André-Quillen cohomology, see \cite[III 2.1.2.3]{ill:71}, cf. \cite[2.5]{ile:19bX}.
By a diagram chase it follows that $\delta$ is surjective and (S2) holds for $\df{}{(A,L')}$. Similarly for (ii).
\end{proof}
%%%%%%%%%%%%%%%%%%%%%%% THEOREM %%%%%%%%%%%%%%%%%%%%%%%
\begin{thm}\label{thm.defgrade}
Consider the map $\sigma_{L'}\co \df{}{(A,N)}\ra\df{}{(A,L')}$ in \textup{Lemma \ref{lem.defmap}.}
\begin{enumerate}[leftmargin=2.4em, label=\textup{(\roman*)}]
\item If\, \(\hm{}{A}{N}{M'}=0\) then \(\sigma_{L'}\) is injective\textup{.}

\item If\, \(\xt{1}{A}{N}{M'}=0\) then \(\sigma_{L'}\) is formally smooth\textup{.}

\item Suppose \(\df{}{(A,N)}\) has a versal element $v=(R\ra \ulset{\tn{v}}{\hspace{-0.1em}\mc{A}},\ulset{\tn{v}}{\mc{N}})$ and \(\xt{1}{A}{N}{M'}=0\)\textup{.} Then $(R\ra \ulset{\tn{v}}{\hspace{-0.1em}\mc{A}},\sigma_{L'}(\hspace{0.06em}\ulset{\tn{v}}{\mc{N}}))$ is a versal element for $\df{}{(A,L')}$ and \(\sigma_{L'}\) is smooth\textup{.}
\end{enumerate}
Analogous statements hold for \(\sigma_{L}\co \df{}{(A,N)}\ra\df{}{(A,L)}\) with \(\xt{1}{A}{N}{M}=0\) in \textup{(i)} and \(\xt{2}{A}{N}{M}=0\) in \textup{(ii-iii).}
\end{thm}
%%%%%%%%%%% EXAMPLE %%%%%%%%%%%
\begin{ex}\label{ex.defgrade}
If \(\grade N\geq n+1\) then $\xt{i}{A}{N}{M}=0$ for all $i\leq n$ and any $M$ in $\MCM_A$.
\end{ex}
\begin{proof}
(i) Suppose \(S\) is an object in \(\He\) and \(({}^{i\hspace{-0.1em}}h\co S\ra \ulset{i}{\hspace{-0.1em}\mc{A}},{}^{i\!}\mc{N})\) are deformations of \((A,N)\) to \(S\) for \(i=1,2\). Assume that the images \(({}^{i\hspace{-0.1em}}h,{}^{i\!}\mc{L}')\) under \(\sigma_{L'}\) are isomorphic, identify \({}^{1\!}h\co S\ra\ulset{1}{\hspace{-0.1em}\mc{A}}\) with \(h=\ulset{2}{\hspace{0.06em}h}\co S\ra\ulset{2}{\hspace{-0.1em}\mc{A}}=\mc{A}\), and let \(\beta\co {}^{1\!}\mc{L}'\ra {}^{2\!}\mc{L}'\) denote the isomorphism. 
Let \(S_{n}=S/\fr{m}_{S}^{n+1}\), \(\mc{A}_{n}=\mc{A}\ot S_{n}\) etc. We construct a tower of isomorphisms \(\{\alpha_{n}\co {}^{1\!}\mc{N}_{n}\cong{}^{2\!}\mc{N}_{n}\}\) which commute with the tower \(\{\beta_{n}\co {}^{1\!}\mc{L}'_{n}\ra {}^{2\!}\mc{L}'_{n}\}\) obtained from \(\beta\) and conclude by Lemma \ref{lem.Aapprox} that the deformations \({}^{1\!}\mc{N}\) and \({}^{2\!}\mc{N}\) are isomorphic. The case \(n=0\) is trivial. Given \(\alpha_{n-1}\) and use it to identify the \({}^{i\!}\mc{N}_{n-1}\) and denote them by \(\mc{N}_{n-1}\). Let \(I=\ker\{S_{n}\ra S_{n-1}\}\). 
The `difference' of the \({}^{i\!}\mc{N}_{n}\) is an element $\gamma$ in \(\xt{1}{A}{N}{N}\ot_{k}I\) by Lemma \ref{lem.omap} which \(\iota_{*}\) maps to \(0\) in \(\xt{1}{A}{N}{L'}\ot_{k}I\). Since \(\iota_{*}\) is injective by assumption, $\gamma=0$.
By Proposition \ref{prop.obsmodule} the \({}^{i\!}\mc{N}_{n}\) are isomorphic by an isomorphism \(\ulset{\star}{\alpha_{n}}\) compatible with \(\alpha_{n-1}\).  Then \(\beta_{n}{}^{1\!}\iota_{n}-{}^{2\!}\iota_{n}\ulset{\star}{\alpha_{n}}\) by the induction hypothesis factors through a \(\delta_{n}\co N\ra L'\ot_{k}I\) which (since \(\hm{}{A}{N}{M'}=0\)) factors through a map \(\eta\co N\ra N\ot I\). Adding the map induced from \(\eta\) to \(\ulset{\star}{\alpha_{n}}\) gives \(\alpha_{n}\) which commutes with $\beta_{n}$.
 
(ii) Let \(S\ra \wbar{S}\) in \(\QA_{k}\) be surjective with kernel \(I\), \(b=(h\co S\ra \mc{A},\mc{L}')\) a deformation of \((A,L')\) to \(S\) and let \(\bar{b}=(\wbar{h}\co \wbar{S}\ra\wbar{\mc{A}},\wbar{\mc{L}}{}')\) denote the base change of \(b\) to \(\wbar{S}\). Suppose there is a deformation \((h^{\star}\co \wbar{S}\ra \mc{A}^{\star},\wbar{\mc{N}})\) of \((A,N)\) which \(\sigma_{L'}\) maps to \(\bar{b}\). As above we can assume that \(h^{\star}=\wbar{h}\). By induction on the length of \(S\) we can assume that \(I\cdot\fr{m}_{S}=0\). By Lemma \ref{lem.omap}, \(\ob(q\co \mc{A}\ra \wbar{\mc{A}},\wbar{\mc{N}})\) maps to \(\ob(q,\wbar{\mc{L}}{}')\) under \(\xt{2}{A}{N}{N}\ot I\ra \xt{2}{A}{L'}{L'}\ot I\) which by the assumption is injective. Since \(\mc{L}'\) lifts \(\wbar{\mc{L}}{}'\) to \(\mc{A}\), \(\ob(q,\wbar{\mc{L}}{}')=0\). By Proposition \ref{prop.obsmodule} there exists a lifting \(\ulset{\star}{\mc{N}}\) of \(\wbar{\mc{N}}\) to \(\mc{A}\). Put \(\ulset{\star}{\mc{L}'}=\sigma_{L'}(\ulset{\star}{\mc{N}})\). The difference of \(\ulset{\star}{\mc{L}'}\) and \(\mc{L}'\) gives a \(\theta\in \xt{1}{A}{L'}{L'}\ot I\). By assumption \(\xt{1}{A}{N}{N}\ot I\) maps surjectively to \(\xt{1}{A}{L'}{L'}\ot I\) and a lifting of \(\theta\) perturbs \(\ulset{\star}{\mc{N}}\) to a lifting \(\mc{N}\) of \(\wbar{\mc{N}}\) with \(\sigma_{L'}(\mc{N})=\mc{L}'\) by Lemma \ref{lem.omap}.

(iii) By Lemma \ref{lem.S2v}, the versality of $v$ implies (S2) for $\df{}{(A,N)}$. Then (S2) follows for $\df{}{(A,L')}$ by Lemma \ref{lem.S2}. Put $\ulset{\tn{v}}{\mc{L}'}=\sigma_{L'}(\hspace{0.06em}\ulset{\tn{v}}{\mc{N}})$ and $v'=(R\ra \ulset{\tn{v}}{\hspace{-0.1em}\mc{A}},\ulset{\tn{v}}{\mc{L}'})$. By (ii), $v'$ is formally versal. To test $v'$ for versality, let \(S\ra S_0\) in \(\He\) be surjective with kernel \(I\) and \(b_0=(h_0\co S_0\ra\mc{A}_0,\mc{L}_0')\) a deformation of \((A,L')\) to \(S_0\) induced from $v'$ by a map $f_0\co R\ra S_0$. Let $b=(h\co S\ra \mc{A},\mc{L}')$ be a lifting of $b_0$ to $S$. Put $S_n=S/I^{n+1}$ and $b_n=b_{S_n}$. As noted by H. van Essen \cite[p. 416]{ess:90}, H. Flenner's \cite[3.2]{fle:81} (where (S2) is needed) implies that a lifting $f\co R\ra S$ of $f_0$ with $f_*v'\cong b$ above $b_0$ exists in the case $I$ is nilpotent; cf. \cite[3.3]{ile:19bX}. This implies that we can find a projective system of maps $\{f_n\co R\ra S_n\}$ and isomorphisms $\{(f_n)_*v'\cong b_n\}$. Let $\hat{f}\co R\ra \limproj S_n=:S_I\hspace{-0.18em}\hat{}\;$  denote the induced map.  
The isomorphism $\limproj \ulset{\tn{v}}{\hspace{-0.1em}\mc{A}}_{S_n}\cong \limproj \mc{A}_{S_n}$ implies that the completions in maximal ideals are isomorphic too; $\ulset{\tn{v}}{\hspace{-0.1em}\mc{A}}\hspace{0.1em}\hat{}\;\cong \mc{A}\hspace{0.1em}\hat{}$\,.

Any \(S\) in \(\He\) is a direct limit of a filtering system of algebraic \(\vL\)-algebras in \(\He\). Since \(\df{}{(A,L')}\) is locally of finite presentation it is sufficient to prove the lifting property for \(S\) algebraic. Since \(\vL\) is excellent, so is \(S\) by \cite[7.8.3]{EGAIV2} and \cite[18.7.6]{EGAIV4}. By Artin's Approximation Theorem \cite[2.6]{art:69} (and \cite[1.3]{pop:86}, \cite{pop:90}) there is an isomorphism $\ulset{\tn{v}}{\hspace{-0.1em}\mc{A}}_{S_I\hspace{-0.18em}\hat{}}\cong \mc{A}_{S_I\hspace{-0.18em}\hat{}}$ over $\ulset{\tn{v}}{\hspace{-0.1em}\mc{A}}_{S_0}\cong \mc{A}_{0}$. By Lemma \ref{lem.Aapprox} there is a corresponding isomorphism of the modules $\ulset{\tn{v}}{\mc{L}'}_{\hspace{-0.24em}S_I\hspace{-0.18em}\hat{}}\cong \mc{L}'_{S_I\hspace{-0.18em}\hat{}}$ compatible with $\ulset{\tn{v}}{\mc{L}'}_{\hspace{-0.24em}S_0}\cong \mc{L}'_{0}$. Hence we have an isomorphism of deformations $f_*v'\cong b_{S_I\hspace{-0.18em}\hat{}}$ above $b_0$. By using Artin's Approximation Theorem \cite[1.12]{art:69} one shows that there is a map $g\co R\ra S$ lifting $f_0$ and an isomorphism of deformations $g_*v'\cong b$ above $b_0$. Smoothness of $\sigma_{L'}$ is equivalent to the versality. The last part is similar.
\end{proof}
An analogous proof gives:
%%%%%%%%%%%%%%%%%%%%%%% THEOREM %%%%%%%%%%%%%%%%%%%%%%%
\begin{thm}\label{thm.defgrade2}
Consider the map \(\sigma_{M}\co \df{}{(A,N)}\lra\df{}{(A,M)}\) in \textup{Lemma \ref{lem.defmap}.}
\begin{enumerate}[leftmargin=2.4em, label=\textup{(\roman*)}]
\item If\, \(\hm{}{A}{L}{N}=0\) then \(\sigma_{M}\) is injective\textup{.}

\item If\, \(\xt{1}{A}{L}{N}=0\) then \(\sigma_{M}\) is formally smooth\textup{.}

\item Suppose \(\df{}{(A,N)}\) has a versal element $(R\ra \ulset{\tn{v}}{\hspace{-0.1em}\mc{A}},\ulset{\tn{v}}{\mc{N}})$ and \(\xt{1}{A}{L}{N}=0\)\textup{.} Then $(R\ra \ulset{\tn{v}}{\hspace{-0.1em}\mc{A}},\sigma_{M}(\hspace{0.06em}\ulset{\tn{v}}{\mc{N}}))$ is a versal element for $\df{}{(A,M)}$ and
\(\sigma_{M}\) is smooth\textup{.}
\end{enumerate}
The analogous statements hold for \(\sigma_{M'}\co \df{}{(A,N)}\ra\df{}{(A,M')}\) with \(\xt{1}{A}{L'}{N}=0\) in \textup{(i)} and \(\xt{2}{A}{L'}{N}=0\) in \textup{(ii-iii).}
\end{thm}
The following two results have very similar proofs to Theorems \ref{thm.defgrade} and \ref{thm.defgrade2}.
%%%%%%%%%%%%%% COROLLARY %%%%%%%%%%%%%%%
\begin{cor}\label{cor.defgrade}
Consider the map \(\sigma^{\mc{A}}_{L'}\co \df{{\mc{A}}}{N}\ra\df{{\mc{A}}}{L'}\) in \textup{Lemma \ref{lem.defmap}.}
\begin{enumerate}[leftmargin=2.4em, label=\textup{(\roman*)}]
\item If\, \(\hm{}{A}{N}{M'}=0\) then \(\sigma^{\mc{A}}_{L'}\) is injective\textup{.}

\item If\, \(\xt{1}{A}{N}{M'}=0\) then \(\sigma^{\mc{A}}_{L'}\) is formally smooth\textup{.}

\item Suppose \(\df{{\mc{A}}}{N}\) has a versal element $(R,\ulset{\tn{v}}{\mc{N}})$ and \(\xt{1}{A}{N}{M'}=0\)\textup{.} Then $(R,\sigma^{\mc{A}}_{L'}(\hspace{0.06em}\ulset{\tn{v}}{\mc{N}}))$ is a versal element for $\df{{\mc{A}}}{L'}$ and
\(\sigma^{\mc{A}}_{L'}\) is smooth\textup{.}
\end{enumerate}
The analogous statements hold for \(\sigma^{\mc{A}}_{L}\co \df{{\mc{A}}}{N}\ra\df{{\mc{A}}}{L}\) with \(\xt{1}{A}{N}{M}=0\) in \textup{(i)} and \(\xt{2}{A}{N}{M}=0\) in \textup{(ii)} and \textup{(iii).}
\end{cor}
%%%%%%%%%%%%%% COROLLARY %%%%%%%%%%%%%%%
\begin{cor}\label{cor.defgrade2}
Consider the map \(\sigma^{\mc{A}}_{M}\co \df{{\mc{A}}}{N}\ra\df{{\mc{A}}}{M}\) in \textup{Lemma \ref{lem.defmap}.}
\begin{enumerate}[leftmargin=2.4em, label=\textup{(\roman*)}]
\item If\, \(\hm{}{A}{L}{N}=0\) then \(\sigma^{\mc{A}}_{M}\) is injective\textup{.}

\item If\, \(\xt{1}{A}{L}{N}=0\) then \(\sigma^{\mc{A}}_{M}\) is formally smooth\textup{.}

\item Suppose \(\df{{\mc{A}}}{N}\) has a versal element $(R\ra \mc{A}_R,\ulset{\tn{v}}{\mc{N}})$ and \(\xt{1}{A}{L}{N}=0\)\textup{.} Then $(R\ra \mc{A}_R,\sigma^{\mc{A}}_{M}(\hspace{0.06em}\ulset{\tn{v}}{\mc{N}}))$ is a versal element for $\df{}{(A,M)}$ and
\(\sigma^{\mc{A}}_{M}\) is smooth\textup{.}
\end{enumerate}
The analogous statements hold for \(\sigma^{\mc{A}}_{M'}\co \df{{\mc{A}}}{N}\ra\df{{\mc{A}}}{M'}\) with \(\xt{1}{A}{L'}{N}=0\) in \textup{(i)} and \(\xt{2}{A}{L'}{N}=0\) in \textup{(ii)} and \textup{(iii).}
\end{cor}
%%%%%%%%%%%%%%%%% PROPOSITION %%%%%%%%%%%%%%%%
\begin{prop}\label{prop.defequiv} 
Put \(Q'=\hm{}{A}{\omega_{A}}{L'}\) and \(Q=\hm{}{A}{\omega_{A}}{L}\). Then:
\begin{enumerate}[leftmargin=2.4em, label=\textup{(\roman*)}]
\item \(Q'\) and \(Q\) have finite projective dimension\textup{.}
\item \(\df{}{(A,L')}\cong\df{}{(A,Q')}\) and \(\df{}{(A,L)}\cong\df{}{(A,Q)}\)\textup{.}
\item There are natural maps 
\begin{equation*}
s\co \df{}{(A,M)}\lra\df{}{(A,M')}\quad \text{and}\quad t\co \df{}{(A,L')}\lra\df{}{(A,L)}
\end{equation*}
commuting with the maps 
\(\sigma_{X}\co \df{}{(A,N)}\ra\df{}{(A,X)}\) for \(X\) equal to \(M\) and \(M'\)\textup{,} and to \(L'\) and \(L\)\textup{,} respectively\textup{.} If \(A\) is a Gorenstein ring\textup{,} then \(s\) is an isomorphism\textup{.}
\end{enumerate}
If $\mc{A}$ is a flat and algebraic $\vL$-algebra with $\mc{A}\ot_{\vL}k\cong A$, the analogous statements hold for the deformation functors \(\df{{\mc{A}}}{X}\)\textup{.}
\end{prop}
\begin{proof}
(i) Applying $\hm{}{A}{\omega_A}{-}$ to a finite $\cat{D}$-resolution of $L'$ gives a finite projective resolution of $Q'$, see \cite[6.10]{ile:12} which also gives (ii).

(iii) There is a short exact sequence \(0\ra M\ra \omega_{A}^{\oplus n}\ra M'\ra 0\) such that the last map is without a common \(\omega_{A}\)-summand, corresponding (through $\omega$-dualisation) to a short exact sequence \(0\la M^{\vee}\la A^{\oplus n}\la (M')^{\vee}\la 0\) where \(n\) is minimal. The map \(s\) is the composition \(\df{}{(A,M)}\cong\df{}{(A,M^{\vee})}\ra\df{}{(A,(M')^{\vee})}\cong\df{}{(A,M')}\) where the first and the last map are given by $\omega$-dualisation. The middle map is given by lifting the surjection $A^{\oplus n}\ra M^{\vee}$ to a free cover of a deformation of $M^{\vee}$ and taking the kernel to get the (minimal) syzygy as a deformation of $(M')^{\vee}$. This is a well-defined map of deformation functors. Then $s$ maps a deformation $\mc{M}\ra M$ to the deformation $(\Syz(\mc{M}^{\vee}))^{\vee}\ra (M')^{\vee\vee}\cong M'$.
If \(A\) is a Gorenstein ring then \(\omega_{A}\cong A\) and $s$ has an inverse \(\df{}{(A,M')}\ra\df{}{(A,M)}\) given by the syzygy map. 

Note that the pushout of \(M\ra \omega_{A}^{\oplus n}\) with \(M\ra N\) gives \(N\ra L'\). Consider the induced short exact sequence \(0\ra L\ra\omega_{A}^{\oplus n}\xra{\;\mu\;} L'\ra 0\). For a deformation \((h,\mc{L}')\) in \(\df{}{(A,L')}\) with structure map \(\lambda'\co \mc{L}'\ra L'\) there is a lifting of \(\mu\) to a map \(\tilde{\mu}\co \omega_{h}^{\oplus n}\ra \mc{L}'\). If \(\mc{L}\) denotes the kernel of \(\tilde{\mu}\) then there is a cocartesian map \(\lambda\co \mc{L}\ra L\) commuting with \(\omega_{h}^{\oplus n}\ra \omega_{A}^{\oplus n}\). By Lemma \ref{lem.defCMapprox2}, \((h,\lambda')\mapsto (h,\lambda)\) gives a well defined map of deformation functors \(t\co \df{}{(A,L')}\ra\df{}{(A,L)}\). 

Given a deformation \((h,\mc{N})\) in \(\df{}{(A,N)}\), let \(0\ra\mc{L}\ra\mc{M}\ra\mc{N}\ra0\) and \(0\ra\mc{N}\ra\mc{L}'\ra\mc{M}'\ra0\) be the minimal sequences in \eqref{eq.ny}.
There is a commutative diagram of short exact sequences with (co)cartesian square (cf.\ \cite{aus/buc:89})
\begin{equation}\label{eq.pullpush}
\xymatrix@C-0pt@R-12pt@H-30pt{
&& 0 \ar[d] & 0 \ar[d] \\
0\ar[r] & \mc{L} \ar[r]\ar@{=}[d] & \mc{M} \ar[r]\ar[d]\ar@{}[dr]|{\Box} & \mc{N} \ar[r]\ar[d] & 0  \\   
0\ar[r] & \mc{L} \ar[r] & \omega_{h}^{\oplus n} \ar[r]\ar[d] & \mc{L}' \ar[r]\ar[d] & 0 \\
&& \mc{M}' \ar@{=}[r]\ar[d] & \mc{M}' \ar[d] \\
&& 0 & 0
}
\end{equation}
where \(\omega_{h}^{\oplus n}\ra \mc{L}'\) is given as above. The stated commutativity of maps of deformation functors follows.
\end{proof}
%%%%%%%%%%%%%% COROLLARY %%%%%%%%%%%%%%%%
\begin{cor}\label{cor.defapprox}
Suppose $A$ has residue field $k$ and \(\dim A\geq 2\)\textup{.} Then there exists finite \(A\)-modules \(L'\) and \(Q'\) with \(\Injdim L'=\dim A=\pdim Q'\) and universal deformations \(\mc{L}'\in\df{A}{L'}(A)\) and \(\mc{Q}'\in\df{A}{Q'}(A)\)\textup{.} 
\end{cor}
\begin{proof}
Let \(h=1\ot\id\co A\ra A\hot_{k}A=\mc{A}\) and \(\mc{N}=A\) be the cyclic \(\mc{A}\)-module defined through the multiplication map. Then \(\mc{A}\ot_{A}k\cong A\) and \(\mc{N}\ot_{A}k\cong k\) and this gives a deformation \(\mc{N}\ra k\) of the residue field of \(A\) which is universal.
If \(L'\) is the minimal \(\hat{\cat{D}}_{A}\)-hull of the residue field \(k\) then \(\mc{L}'=\sigma_{L'}(\mc{N})\in \df{{\mc{A}}}{L'}(A)\) is universal by Corollary \ref{cor.defgrade}. If \(Q'=\hm{}{A}{\omega_{A}}{L'}\) then \(\hm{}{\mc{A}}{\omega_{\mc{A}}}{\mc{L}'}\in\df{A}{Q'}(A)\) is universal by Proposition \ref{prop.defequiv}.
\end{proof}
%%%%%%%%%%%%%% COROLLARY %%%%%%%%%%%%%%%%%
\begin{cor}\label{cor.depth}
Put \(X=\Spec A\)\textup{.} Let \(Z\) be a closed subscheme of \(X\) such that the complement \(U\) is contained in the regular locus\textup{.} Assume \(\tilde{N}_{\vert U}\) is locally free\textup{,} \(\depth_{Z}N\geq 2\) and \(\cH^{2}_{Z}(\hm{}{A}{L}{N})=0\)\textup{.} Then \(\xt{1}{A}{L}{N}=0\) and so
\begin{equation*}
\sigma_{M}\co \df{}{(A,N)}\lra\df{}{(A,M)}\quad\text{and}\quad \sigma_{M}^{{\mc{A}}}\co \df{{\mc{A}}}{N}\lra\df{{\mc{A}}}{M}\quad\text{are formally smooth\textup{.}}
\end{equation*}
\end{cor}
\begin{proof}
We show that \(\xt{1}{A}{L}{N}=0\) and apply Theorem \ref{thm.defgrade2} and Corollary \ref{cor.defgrade2}.
By Theorem 1.6 in \cite[Expos{\'e} VI]{SGA2} there is a cohomological spectral sequence
\begin{equation}
\cE^{p,q}_{2}=\xt{q}{A}{L}{\cH^{p}_{Z}(N)}\,\Ra\,\xt{p+q}{Z}{X;L}{N}\,.
\end{equation}
Since \(\cH^{i}_{Z}(N)=0\) for \(i=0,1\) the restriction map \(\xt{1}{A}{L}{N}\ra\xt{1}{U}{X;L}{N}\) in the long exact sequence is injective. Since \(U\) is contained in the regular locus, \(\tilde{M}_{\vert U}\) and hence \(\tilde{L}_{\vert U}\) are locally free. It follows that \(\xt{1}{U}{X;L}{N}\) is isomorphic to 
\begin{equation}
\xt{1}{\Q_{U}}{\tilde{L}_{\vert U}}{\tilde{N}_{\vert U}}\cong\cH^{1}(U,\shm{}{\Q_{\hspace{-0.1em}X}\hspace{-0.2em}}{\tilde{L}}{\tilde{N}}) \cong \cH^{2}_{Z}(\hm{}{A}{L}{N})
\end{equation}
which is zero by assumption.
\end{proof}
%%%%%%%%% EXAMPLE %%%%%%%%%%%
\begin{ex}\label{ex.depth}
The condition \(\cH^{2}_{Z}(\hm{}{A}{L}{N})=0\) is implied by \(\tilde{N}_{\vert U}=0\) and also by \(\depth_{Z}(\hm{}{A}{L}{N})\geq 3\).
\end{ex}
The following result extends A.\ Ishii's \cite[3.2]{ish:00} to deformations of the pair.
%%%%%%%%%%%%%% PROPOSITION %%%%%%%%%%%%%%%
\begin{prop}\label{prop.Gor}
Assume \(A\) is Gorenstein\textup{.} If \(\depth N=\dim A -1\) then 
\begin{equation*}
\sigma_{M}\co \df{}{(A,N)}\lra\df{}{(A,M)}\quad\text{and}\quad \sigma_{M}^{{\mc{A}}}\co \df{{\mc{A}}}{N}\lra\df{{\mc{A}}}{M}\quad\text{are smooth\textup{.}}
\end{equation*}
\end{prop}
\begin{proof}
Let \(S_2\ra S_1\) be a surjection in \(\He\) and \((h_2\co S_2\ra \mc{B}_2,\mc{M}_{2})\) an element in \(\df{}{(A,M)}(S_2)\) which maps to \((h_1\co S_1\ra \mc{B}_1,\mc{M}_1)\) in \(\df{}{(A,M)}(S_1)\). Suppose \(\sigma_{M}\) maps \((h',\mc{N}_1)\) in \(\df{}{(A,N)}(S)\) to \((h_1,\mc{M}_1)\). 
By the depth lemma, $\depth L=\depth N+1=\dim A$, so $L$ is a MCM-module of finite injective dimension, hence \(L\cong A^{\oplus r}\) for some \(r\) since $A$ is Gorenstein (see \cite[3.7]{aus/buc:89} for a more general statement). 
We can assume that \(h'=h_1\) and that the minimal \(\MCM\)-approximation of \(\mc{N}\) is \(0\ra\mc{L}_1\xra{\rho_1}\mc{M}_1\ra\mc{N}_1\ra0\) where \(\mc{L}_1\cong \mc{B}_1^{\oplus r}\). Put \(\mc{L}_{2}:=\mc{B}_2^{\oplus r}\) and choose a lifting \(\rho_{2}\co \mc{L}_{2}\ra\mc{M}_{2}\) of \(\rho_1\). Put \(\mc{N}_{2}:=\coker\rho_{2}\) with its natural map to \(\mc{N}_1\). Then \(\mc{N}_{2}\) is \(S_2\)-flat (\(\rho_{2}\ot S_1=\rho_1\)) and \(\sigma_{M}(h_2,\mc{N}_{2})=(h_2,\mc{M}_{2})\). 
\end{proof}
%%%%%%%%%%% REMARK %%%%%%%%%%%
\begin{rem}
If \(A\) is a Gorenstein domain and \(M\) is an MCM \(A\)-module there is a short exact sequence \(0\ra A^{\oplus r}\ra M\ra N\ra 0\) with \(N\) a codimension one Cohen-Macaulay module; cf.\ \cite[1.4.3]{bru/her:98}. This sequence is an \(\MCM_{A}\)-approximation and Proposition \ref{prop.Gor} applies. However, it is not always possible to continue this reduction. Assume \(A\) is a normal Gorenstein complete local ring. Then all MCM \(A\)-modules are \(\MCM_{A}\)-approximations of codimension \(2\) Cohen-Macaulay modules up to stable isomorphism if and only if \(A\) is a unique factorisation domain; see \cite{yos/iso:00, kat:07}.

Let \(A\) be a Gorenstein normal domain of dimension \(2\) and \(0\ra A^{\oplus r{-}1}\ra M\ra I\ra 0\) the minimal MCM approximation of a torsion-free rank \(1\) module \(I\). Let \(U\) denote the regular locus in \(X=\Spec A\). If \(\mc{A}=A\hot_{k}S\) for \(S\) in \(\kH_k\) there is a natural section \(A\ra \mc{A}\). Let \(U_{\!\mc{A}}\) denote \(U{\times}_{X}\Spec \mc{A}\). Consider the subfunctor \(\df{A,\wedge}{M}\sbeq\df{A}{M}\) of deformations \(\mc{M}\) such that \(\wedge^{r}\mc{M}_{\vert U_{\!\mc{A}}}\cong\Q_{U_{\!\mc{A}}}\). Note that \(\cH^{0}(U,\wedge^{r}M)\) is isomorphic to the MCM \(A\)-module \(\wbar{I}:=\cH^{0}(U,I)\). Proposition \ref{prop.Gor} implies that the resulting map from the (local) functor of quotients \(\quot{\wbar{I}}{I\sbeq \wbar{I}}\ra \df{A,\wedge}{M}\) is smooth; cf.\ \cite[3.2]{ish:00}. In particular, if \(E_{A}\) is the fundamental module (see \eqref{eq.funda} below) and \(A/\fr{m}_{A}\cong k\) then \(h_A\cong\quot{A}{\fr{m}_{A}\sbeq A}\cong\df{A}{k}\) gives a versal family for \(\df{A,\wedge}{E_{A}}\) by the MCM approximation in \cite[7.4]{ile:12}; see \cite[3.4]{ish:00}.
\end{rem}
%%%%%%%%%%% EXAMPLE %%%%%%%%%%%%
\begin{ex}
Assume \(A/\fr{m}_{A}\cong k\) and let \(M\) denote the minimal MCM approximation of \(k\). It is given as \(M\cong \hm{}{A}{\syz{A}{d}(k^{\vee})}{\omega_{A}}\) where \(d=\dim A\); cf.\ \cite[5.6]{ile:12}. One has \(k^{\vee}=\xt{d}{A}{k}{\omega_{A}}\cong k\). 
 We apply \(\hm{}{A}{-}{\omega_{A}}\) to the short exact sequence \(0\ra \syz{A}{}(\fr{m}_{A})\ra A^{\oplus\beta_{1}}\xra{(\ul{x})}\fr{m}_{A}\ra 0\). Assume \(\dim A = 2\). Since \(\xt{1}{A}{\fr{m}_{A}}{\omega_{A}}\cong k\) we obtain the MCM approximation of \(k\) from the exact sequence
\begin{equation}
0\ra\omega_{A}\xra{(\ul{x})^{\text{tr}}}\omega_{A}^{\oplus\beta_{1}}\lra M\lra k\ra 0\,.
\end{equation}
In particular \(\rk(M)=\beta_{1}-1\). Put $\mu(M)=\dim M/\fr{m}_AM$ and let \(t(A)\) denote the type of \(A\); cf. \cite[1.2.15]{bru/her:98}. Then \(\mu(M)=t(A)\cdot\beta_{1}+1\); cf. \cite[3.3.11]{bru/her:98}. 

If $k=\bar{k}$ and \(A=A(m)=k[u^{m},u^{m-1}v,\dots,v^{m}]^{\text{h}}\), the vertex of the cone over the rational normal curve of degree \(m\), the indecomposable MCM \(A\)-modules are \(M_{i}=(u^{i},u^{i-1}v,\dots,v^{i})\) for \(i=0,\dots,m{-}1\) (an argument independent of characteristic is given in \cite[Lemma 1]{gus/ile:04b}). In particular $\omega_A=M_{m-2}$ (cf. \cite[p. 616]{gus/ile:04b}) so that $\mu(M)=m^2$ and from this \(M=M_{m-1}^{\oplus m}\) follows. Then 
\begin{equation}
\dim_{k}\df{A}{M}(k[\vare])=\dim_{k}\xt{1}{A}{M}{M}= (m-1)\cdot m^{2}
\end{equation}
by applying a result of Ishii; cf. equation (12) and calculations on p. 616 in \cite{gus/ile:04b}. Moreover, applying $\hm{}{A}{-}{k}$ to $0\ra\fr{m}\ra A\ra k\ra 0$ gives \(\dim_{k}\df{A}{k}(k[\vare])=\beta_1=m+1\). Even in the Gorenstein case (\(m=2\)) the tangent map is not surjective and so Proposition \ref{prop.Gor} cannot in general be extended to \(\depth N=\dim A -2\). See \cite{gus/ile:04b} for a detailed description of the strata of the reduced versal deformation space of \(M\) defined by Ishii in \cite{ish:00}. 

If \(\dim A=2\) the \(\MCM_{A}\)-approximation of \(\fr{m}_{A}\) is a short exact sequence 
\begin{equation}\label{eq.funda}
0\ra\omega_{A}\lra E_{A}\lra \fr{m}_{A}\ra0
\end{equation}
where \(E_{A}\) is called the \emph{fundamental module}; cf.\ \cite[7.1.1]{ile:12}. Applying \(\hm{}{A}{k}{-}\) to \(0\ra\fr{m}_{A}\ra A\ra k\ra0\) gives an exact sequence
\begin{equation}
0\ra\xt{1}{A}{k}{k}\lra\df{A}{\fr{m}_{A}}(k[\vare])\lra k^{\oplus t(A)}\lra \xt{2}{A}{k}{k}
\end{equation}
since \(\xt{1}{A}{\fr{m}_{A}}{\fr{m}_{A}}\cong\xt{2}{A}{k}{\fr{m}_{A}}\) and \(\dim A = 2\). If \(A=A(m)\) then \(E_{A}\) is isomorphic to \(M_{m-1}^{\oplus 2}\) with \(\dim_{k}\df{A}{E_{A}}(k[\vare])=4(m-1)\). Hence the conclusion in Proposition \ref{prop.Gor} cannot hold in the non-Gorenstein case \(m>2\).
\end{ex}
%%%%%%%%%%%%%%%%%%%%%%%%%%%%%%%%%%%%%%%%%%%%%%
%%%%%%%%%%%%%%%%%%%% SECTION %%%%%%%%%%%%%%%%%%%%%
%%%%%%%%%%%%%%%%%%%%%%%%%%%%%%%%%%%%%%%%%%%%%%
\section{Deforming maximal Cohen-Macaulay approximations of Cohen-Macaulay modules}\label{sec.defCM}
Several definitions and results are given to prepare the statement of Theorem \ref{thm.defMCM} and then to prove it.
%%%%%%%%% DEFINITION %%%%%%%%%%%%%
\begin{defn}\label{def.obsth}
A functor \(F\co \QA_{k}\ra \Sets\) has an \emph{obstruction theory} if there is a $k$-linear functor \(\cH_{F}^{2}\co\Mod_{k}\ra\Mod_{k}\) and for each small surjection \(p\co R\ra S\) in \(\QA_{k}\) (i.e. with kernel \(I\) such that \(\fr{m}_{R}{\cdot} I=0\)) and each \(a\in F(S)\) there is an element \(\mr{o}(p,a)\in \cH_{F}^{2}(I)\) which is zero if and only if there exists a \(b\in F(R)\) mapping to \(a\). The obstruction should be functorial with respect to such lifting situations.
Cf. \cite[2.6]{art:74}. 
\end{defn} 
%%%%%%%% EXAMPLE %%%%%%%%%
\begin{ex}\label{ex.ot}
Consider the functor $F=\df{\mc{A}}{N}\co \He\ra \Sets$ defined in Section \ref{sbsec.def} where $N$ is an $A=\mc{A}\ot_\vL k$-module. For a small surjection \(p\co R\ra S\) with kernel $I$, let $J$ denote the kernel of the induced $q\co\mc{A}_R\ra \mc{A}_S$. If $\mc{N}$ is a deformation of $N$ to $S$, there is an obstruction element $\ob(q,\mc{N})$ in $\xt{2}{\mc{A}_S}{\mc{N}}{\mc{N}\ot J}\cong \xt{2}{A}{N}{N}\ot_k I$ which is natural for the lifting situation by Proposition \ref{prop.obsmodule} and Lemma \ref{lem.obsdef}. Then $\cH^2_F(-):=\xt{2}{A}{N}{N}\ot_k(-)$ with the obstruction $\mr{o}(p,\mc{N}):=\ob(q,\mc{N})$ gives an obstruction theory for $\df{\mc{A}}{N}$.
\end{ex}
%%%%%%%%%%%%%%%% LEMMA %%%%%%%%%%%%%%%%
\begin{lem}\label{lem.ind}
Given a map $f\co R\ra S$ in $\He$ with both rings being algebraic over $\vL$ \textup{(}or complete\textup{)} such that the induced map 
$$
\fr{m}_{R}/(\fr{m}_{R}^{2}\!+\im\fr{m}_{\vL}{\cdot}R)\lra \fr{m}_{S}/(\fr{m}_{S}^{2}\!+\im\fr{m}_{\vL}{\cdot}S)
$$ is surjective\textup{.} Then $f$ is a surjection\textup{.}
\end{lem}
\begin{proof}
Let $t_{R/\vL}$ denote the relative Zariski tangent space $[\fr{m}_{R}/(\fr{m}_{R}^{2}\!+\im\fr{m}_{\vL}{\cdot}R)]^*$. There is a $\vL$-algebra map $f_{\tn{lc}}\co R_{\tn{lc}}\ra S_{\tn{lc}}$ which is the Zariski localisation in $k$-points of a map of finite type $\vL$-algebras such that the henselisation of $f_{\tn{lc}}$ is $f$. The induced map $t^*_{R_{\tn{lc}}/\vL}\ra t^*_{R/\vL}$ is an isomorphism, and likewise for $S$. Then  $t^*_{R_{\tn{lc}}/\vL}\ra t^*_{S_{\tn{lc}}/\vL}$ surjective implies $\fr{m}_{R_{\tn{lc}}}/(\fr{m}_{R_{\tn{lc}}})^2 \ra \fr{m}_{S_{\tn{lc}}}/\fr{m}_{S_{\tn{lc}}}^2$ surjective;  
cf. \cite[\href{http://stacks.math.columbia.edu/tag/06GB}{Tag 06GB}]{SP}. Then $\im\fr{m}_{R_{\tn{lc}}}\!{\cdot}\, S_{\tn{lc}}=\fr{m}_{S_{\tn{lc}}}$ by Nakayama's lemma. In particular, $S_{\tn{lc}}$ is the Zariski localisation of a finite $R_{\tn{lc}}$-algebra $S_{1}$ by \cite[\href{http://stacks.math.columbia.edu/tag/052V}{Tag 052V}]{SP}.  Since $\hat{R}_{\tn{lc}}\ra \hat{S}_{\tn{lc}}\cong \hat{S}_{1}$ is surjective by \cite[\href{http://stacks.math.columbia.edu/tag/0315}{Tag 00M9}]{SP}, $R_{\tn{lc}}\ra S_{1}$ is surjective by faithfully flatness of completion. Since henselisation preserves surjections $f$ is surjective.
\end{proof}
%%%%%%%%%%%%%%
\begin{ex}
Suppose $F\co\He\ra\Sets$ is a functor with versal elements in $F(R)$ and $F(S)$ such that the induced maps $h_R(k[\vare])\ra F(k[\vare])\la h_S(k[\vare])$ are bijective. 
Then $R\cong S$. Indeed, by versality there are maps $f\co R\ra S$ and $g\co S\ra R$ which are surjections by Lemma \ref{lem.ind}. Then $gf$ is an automorphism since $R$ is noetherian.
\end{ex}
Put $t_{F\hspace{-0.1em}/\hspace{-0.1em}\vL}=F(k[\vare])$. In the case $k_0\ra k$ is a separable field extension, we will call a base ring $R$ of a (formally) versal (formal) element in $F$ for \emph{minimal} if the induced map $h_R(k[\vare])\ra t_{F\hspace{-0.1em}/\hspace{-0.1em}\vL}$ is bijective; cf. \cite[\href{https://stacks.math.columbia.edu/tag/06IL}{Tag 06IL}]{SP}.
%%%%%%%%%%%%%%%% LEMMA %%%%%%%%%%%%%%%%
\begin{lem}\label{lem.res}%3
Suppose $k_0\ra k$ is a separable field extension and \(\phi\co F\ra G\) is a map of set\textup{-}valued functors on $\QA_{k}$ which have minimal formally versal formal families with base rings \(R^{F}\) and \(R^{G}\) which are algebraic over $\vL$ \textup{(}or complete\textup{).} Put \(V=\ker\{t_{G\hspace{-0.1em}/\hspace{-0.1em}\vL}^{*}\ra t_{F\hspace{-0.1em}/\hspace{-0.1em}\vL}^{*}\}\)\textup{.} Assume\textup{:}
\begin{enumerate}[leftmargin=2.4em, label=\textup{(\roman*)}]
\item The map \(t_{F\hspace{-0.1em}/\hspace{-0.1em}\vL}\ra t_{G\hspace{-0.1em}/\hspace{-0.1em}\vL}\) is injective\textup{.}
\item There are obstruction theories for \(F\) and \(G\) such that \(\mr{o}_G(p,\phi_{S}(\zeta))=0\) implies \(\mr{o}_F(p,\zeta)=0\) for any small surjection \(p\co R\ra S\) in $\QA_{k}$ and element \(\zeta \in F(S)\)\textup{.}
\end{enumerate}
Then every \(f\co R^{G}\ra R^{F}\) in \(\He\) lifting \(\phi\) is surjective and the ideal \(\ker\hspace{-0.1em}f\) is generated by a lifting of a \(k\)\textup{-}basis for \(V\)\textup{.} In particular \(\ker\hspace{-0.1em} f\) is generated by \textup{`}linear forms\textup{'} modulo \(\im\fr{m}_{\vL}{\cdot}R^{G}\)\textup{.}
\end{lem}
\begin{proof}
The Jacobi-Zariski-sequence of an object $\vL\ra R\ra k$ in $\He$ gives the exact sequence (cf. \cite[\href{http://stacks.math.columbia.edu/tag/06S9}{Tag 06S9}]{SP})
\begin{equation}\label{eq.JZ}
\xymatrix@C+3pt@R-6pt@H-30pt{
(\fr{m}_\vL/\fr{m}_\vL^2)\ot_{k_0}k\ar[r] & \fr{m}_R/\fr{m}_R^2 \ar[r]^(0.43){d} & \Omega_{R/\vL}\ot_Rk\ar[r] & \Omega_{k/\vL}\ar[r] & 0
} 
\end{equation}
where $\Omega_{k/\vL}\cong \Omega_{k/k_0}$ which equals $0$ by separability. 
Then 
\begin{equation}\label{eq.cotang}
\Omega_{R/\vL}{\ot}_Rk\cong \fr{m}_{R}/(\fr{m}_{R}^{2}\!+\im\fr{m}_{\vL}{\cdot}R)=t_{R/\vL}^*
\end{equation}
Hence 
\begin{equation}
h_R(k[\vare])\cong\der{\vL}{R}{k}\cong\hm{}{k}{\Omega_{R/\vL}{\ot}_Rk}{k}\cong t_{R/\vL}.
\end{equation}
Then the surjective map \(\phi(k[\vare])^*\co t_{G\hspace{-0.1em}/\hspace{-0.1em}\vL}^{*}\ra t_{F\hspace{-0.1em}/\hspace{-0.1em}\vL}^{*}\) by minimality (cf.  \cite[\href{http://stacks.math.columbia.edu/tag/06IL}{Tag 06IL}]{SP}) is canonically isomorphic to the map
$
t^*_{R^G/\vL}\ra t^*_{R^F/\vL}
$
induced by $f$ so \(f\) is surjective by Lemma \ref{lem.ind}. Moreover, \(\ker\hspace{-0.1em}f\) maps surjectively to \(V\) inducing the natural surjective \(k\)-linear map 
\begin{equation}
g\co(\ker\hspace{-0.1em} f)/(\fr{m}_{R^{G}}\hspace{-0.2em}\cdot\ker f+\im\fr{m}_{\vL}{\cdot}R^{G}\cap\ker\hspace{-0.1em}f)\lra V\,.
\end{equation}
Lift a \(k\)-basis for \(V\) to elements in \(\ker\hspace{-0.1em}f\) and let \(J\) be the ideal in \(R^{G}\) generated by these elements. Then \(g\) is an isomorphism if and only if \(J=\ker\hspace{-0.1em}f\). Put \(\wbar{R}=R^{G}\hspace{-0.1em}/J\), \(\wbar{R}_{n}=\wbar{R}/(\fr{m}_{\wbar{R}}^{n+1}\!+\im\fr{m}_{\vL}\cdot\fr{m}_{\wbar{R}}^{n-1})\). Let \(\wbar{\zeta}_{n}\in G(\wbar{R}_{n})\) for \(n=1,2,\dots\) denote the images of a formal versal family \((\zeta_{n})\) for \(G\). Similarly, put \(R^{F}_{n}:=\wbar{R}/(\fr{m}_{R^{F}}^{n+1}\!+\im\fr{m}_{\vL}\cdot\fr{m}_{R^{F}}^{n-1})\) and let \((\xi_{n})\), \(\xi_{n}\in F(R^{F}_{n})\), denote a formal versal family for \(F\). We prove that the maps \(\wbar{R}_{n}\ra R^{F}_{n}\) are isomorphisms by induction on \(n\). Surjectivity and isomorphic completions imply that \(\wbar{R}\ra R^{F}\) is an isomorphism also in the algebraic case. Put \(K_{1}=\ker\{R^{G}_{1}\ra \wbar{R}_{1}\}\). Then \(K_{1}\) is contained in \(V\), but since \(J\ra V\) is surjective and factors through \(K_{1}\) we have \(K_{1}=V\). This is equivalent to \(t_{\wbar{R}\hspace{-0.0em}/\hspace{-0.1em}\vL}^{*}\cong t_{F\hspace{-0.1em}/\hspace{-0.1em}\vL}^{*}\) and to \(\wbar{R}_{1}\cong R^{F}_{1}\). Let \(\wbar{f}_{n}\co \wbar{R}_{n}\ra R^{F}_{n}\) be the map induced from \(f\). Assume \(\wbar{f}_{n-1}:\wbar{R}_{n-1}\cong R^{F}_{n-1}\). Then \(\wbar{R}_{n}\ra R^{F}_{n-1}\) is a small surjection in $\QA_{k}$ and by (ii) there is an element \(\eta \in F(\wbar{R}_{n})\) lifting \(\xi_{n-1}\). By formal versality there is a map \(h'\co R^{F}_{n}\ra \wbar{R}_{n}\) above \(R^{F}_{n-1}\) such that \(F(h')(\xi_{n})=\eta\). Then \(\wbar{f}_{n}h'\) is an automorphism of \(R^{F}_{n}\) lifting the identity on \(R^{F}_{n-1}\). Precomposing \(h'\) with the inverse of this automorphism gives a section \(h\) to \(\wbar{f}_{n}\). Then \(h\) is surjective too and \(\wbar{f}_{n}\) is an isomorphism.
\end{proof}
Let \(h\co S\ra \mc{A}\) be a flat and local map of noetherian rings. An \emph{\(h\)-sequence} is a sequence \(\mc{J}=(f_{1},\dots,f_{n})\) in \(\mc{A}\) such that the image \({J}\) in \(A=\mc{A}\ot_{S}S/\fr{m}_{S}\) is an \(A\)-sequence. By \cite[2.5]{ile:14} \(\mc{J}\) is an \(h\)-sequence if and only if \(\mc{J}\) is an \(\mc{A}\)-sequence and \(\mc{A}/\hspace{-0.2em}\mc{J}\) is \(S\)-flat. I.e.\ \(\mc{J}\) is a transversally \(\mc{A}\)-regular sequence relative to \(S\) as defined in \cite[19.2.1]{EGAIV4}. 
%%%%%%%%%%%%%%% THEOREM %%%%%%%%%%%%%%%%
\begin{thm}\label{thm.defMCM}
Suppose $k_0\ra k$ is a separable field extension.
Let \(h\co \vL\ra {\mc{A}}\) denote the henselisation of a flat and finite type ring map at a maximal ideal \textup{(cf. Section \ref{subsec.loc}).} Assume \(A={\mc{A}}\ot_{\hspace{-0.12em}\vL}k\) is Cohen-Macaulay and \(\mc{J}=(f_{1},\dots,f_{n})\) is an \(h\)\textup{-}sequence\textup{.} Put \({\mc{B}}=\mc{A}/\!\mc{J}\)\textup{,} \(B=\mc{B}\ot_{\hspace{-0.12em}\vL}k\) and let \(J\) be the image of \(\mc{J}\) in \(A\)\textup{.}
Let \(N\) be a maximal Cohen-Macaulay \(B\)-module and 
\begin{equation*}
0\ra L\lra M\lra N\ra 0
\end{equation*}
a minimal \(\MCM_{A}\)-approximation of \(N\)\textup{.} 

Suppose $\df{\mc{B}}{N}$ and $\df{\mc{A}}{M}$ have formally versal formal  families \textup{(}versal families\textup{)} for minimal base rings \(R^{N}\) and \(R^{M}\) which are complete \textup{(}respectively algebraic over $\vL$\textup{).}
If\, \(\ob{}{}(A/J^{2}\ra B,N)=0\) then 
\begin{equation*}
R^{N}\cong R^{M}\hspace{-0.2em}/J
\end{equation*}
where \(J\) is generated by elements lifting a \(k\)\textup{-}basis of the kernel of the map of dual Zariski vector spaces \textup{(cf. Lemma \ref{lem.res}).} 
In particular \(J\) is generated by \textup{`}linear forms\textup{'} modulo \(\im\fr{m}_{\vL}{\cdot} R^{M}\)\textup{.}
\end{thm}
\begin{ex}
The existence of a splitting of \(q\co A/{J}^{2}\ra B\) implies that \(\ob{}{}(q,N)=0\)
for all \(B\)-modules \(N\) since \(A/{J}^{2}\ot_{B}N\) gives a lifting of \(N\) to \(A/{J}^{2}\).
\end{ex}
Let \(\cat{C}\) be a category. Then \(\Arr\cat{C}\) denotes the category with objects being arrows in \(\cat{C}\) and arrows being commutative diagrams of arrows in \(\cat{C}\). An endo-functor \(F\) on \(\cat{C}\) induces an endo-functor \(\Arr F\) on \(\Arr\cat{C}\). Let \(B\) be a noetherian local ring and \(\cat{P}_{\!B}\) the additive subcategory of projective modules in \(\cat{mod}_{B}\). Let \(\uhm{}{B}{N}{M}\) denote the homomorphisms from \(N\) to \(M\) in the quotient category \(\umod_{B}=\cat{mod}_{B}/\cat{P}_{\!B}\) i.e.\ \(B\)-homomorphisms modulo the ones factoring through an object in \(\cat{P}_{\!B}\). For each \(N\) in \(\cat{mod}_{B}\) we fix a minimal \(B\)-free resolution and use it to define the syzygy modules of \(N\). For each \(i\) the association \(N\mapsto\syz{B}{i}N\) induces an endo-functor on \(\umod_{B}\) defined by A.\ Heller \cite{hel:60}. Define \(\uxt{i}{B}{N}{M}\) as \(\uhm{}{B}{\syz{B}{i}N}{M}\) which turns out to be isomorphic to \(\xt{i}{B}{N}{M}\) for all \(i>0\).
%%%%%%%%%%%%% LEMMA %%%%%%%%%%%%
\begin{lem}\label{lem.obssplit}
Let \(A\) be a noetherian local ring and \(J=(f_{1},\dots,f_{n})\) a regular sequence\textup{.} Put \(B=A/J\) and suppose \(N\) and $N_j$ \textup{(}$j=1,2$\textup{)} are finite \(B\)-modules\textup{.} Let \(\wbar{M}_{\hspace{-0.15em}N}\) denote \(B\ot_{A}\syz{A}{n}N\)\textup{.}
\begin{enumerate}[leftmargin=2.4em, label=\textup{(\roman*)}]
\item There is an injective map \(u_{N}\co N\ra \wbar{M}_{\hspace{-0.15em}N}
\) of \(B\)-modules which induces a functor \(u\co \ul{\cat{mod}}_{B}\ra\Arr\ul{\cat{mod}}_{B}\)\textup{.}
\item The functor \(u\) commutes with the \(B\)-syzygy functor\textup{:} 
\begin{equation*}
\Arr\syz{B}{i}(u_{N})=u_{\syz{B}{i}N}
\end{equation*}
\item Put $\wbar{M}_j=\wbar{M}_{\hspace{-0.15em}N_{j}}$\textup{.} The endo-functor \(B\ot_{A}\syz{A}{n}(-)\) induces a natural map 
\begin{equation*}
\uxt{i}{B}{N_{1}}{N_{2}}\ra\uxt{i}{B}{\wbar{M}_{1}}{\wbar{M}_{2}}
\end{equation*}
which makes the following diagram commutative for all \(i\)\textup{:}
\begin{equation*}
\xymatrix@C-12pt@R-12pt@H-0pt{
\uxt{i}{B}{N_{1}}{N_{2}}\ar[rr]\ar[dr]_(0.4){(u_{N_{2}})_{*}} && \uxt{i}{B}{\wbar{M}_{1}}{\wbar{M}_{2}}\ar[dl]^(0.4){(u_{N_{1}})^{*}}\\
& \uxt{i}{B}{N_{1}}{\wbar{M}_{2}} &
}
\end{equation*}
\item The inclusion \(u_{N}\co N\hra B\ot_{A}\syz{A}{n}N\) splits \(\Llra\ob(A/J^{2}\ra B,N)=0\).
\end{enumerate}
\end{lem}
%%%%%%%%%%%%% REMARK %%%%%%%%%%%%
\begin{rem}
Lemma \ref{lem.obssplit} (iv) strengthens \cite[3.6]{aus/din/sol:93} (in the commutative case). 
\end{rem}
\begin{proof}
(i) Suppose \(F_{*}\ra N\) is the fixed minimal \(A\)-free resolution of \(N\). Tensoring the short exact sequence \(0\ra\syz{A}{n}N\xra{\iota} F_{n-1}\ra\syz{A}{n{-}1}N\ra 0\) with \(B\) gives the exact sequence 
\begin{equation}
0\ra\tor{A}{1}{B}{\syz{A}{n{-}1}N}\ra \wbar{M}_{\hspace{-0.15em}N}\ra B\ot_{\hspace{-0.15em}A}F_{n-1}\ra B\ot\syz{A}{n{-}1}N\ra 0\,.
\end{equation}
We have \(\tor{A}{1}{B}{\syz{A}{n{-}1}N}\cong \tor{A}{n}{B}{N}\cong N\) (use the Koszul complex $K(\ul{f})$ to resolve $B$). Let \(u_{N}\) be the composition \(N\cong\ker(B\ot_{\hspace{-0.15em}A}\iota)\sbeq \wbar{M}_{\hspace{-0.15em}N}\). Then \(N\mapsto u_{N}\) gives a functor of quotient categories.

(ii) Let \(p\co Q\ra N\) be the minimal \(B\)-free cover and \(P_{*}\ra \syz{B}{}N\) the minimal \(A\)-free resolution of the \(B\)-syzygy \(\ker(p)\). Then there is an \(A\)-free resolution \(H_{*}\ra Q\) which is an extension of \(F_{*}\) by \(P_{*}\). Since \(\syz{A}{n}B\cong A\), tensoring the short exact sequence of \(A\)-free resolutions \(0\ra P_{*}\ra H_{*}\ra F_{*}\ra0\) by \(B\) we obtain by (i) a commutative diagram with exact rows 
\begin{equation}
\xymatrix@C-4pt@R-6pt@H-0pt@M4pt{
0\ar[r] & \syz{B}{}N\ar[d]_(0.45){u_{\Syz\hspace{-0.1em}N}}\ar[r] & Q\ar[d]\ar[r] & N\ar[d]_(0.4){u_{N}}\ar[r] & 0 \\
0\ar[r] & B\,\ot\,\syz{A}{n}(\syz{B}{}N)\ar[r] & B^{r}{\oplus}\,Q\ar[r] & B\,\ot\,\syz{A}{n}N\ar[r] & 0
}
\end{equation}
which proves the claim.

(iii) By (ii) it is enough to prove this for \(i=0\). The case \(i=0\) follows from the functoriality in (i).

(iv,\(\La\)) For the case \(n=1\) see the proof of \cite[3.2]{aus/din/sol:93}. Assume \(n\geq 2\). We follow the proof of \cite[3.6]{aus/din/sol:93} closely. Let \(A_{1}=A/(f_{1})\). Then \(F_{*}^{(1)}=A_{1}\ot F_{*\geq 1}[1]\) gives a minimal \(A_{1}\)-free resolution of \(A_{1}\ot\syz{A}{}N\). We have \(\ob(A/J^{2}\ra B,N)=0\Ra \ob(A/(f_{1})^{2}\ra A_{1},N)=0\) and hence \(N\) is a direct summand of \(A_{1}\ot\syz{A}{}N\). Let \(G_{*}\ra N\) be a minimal \(A_{1}\)-free resolution of \(N\). Then \(G_{*}\) is a direct summand of \(F_{*}^{(1)}\) and hence \(\syz{A_{1}}{n{-}1}N\) is a direct summand of \(\syz{A_{1}}{n{-}1}(A_{1}\ot\syz{A}{}N)=A_{1}\ot \syz{A}{n}N\). Tensoring this situation with \(B\) (and let $\wbar{F}=B\ot F$) gives a commutative diagram:
\begin{equation}\label{eq.split}
\xymatrix@C-4pt@R-6pt@H-30pt{
N\ar[r]_(0.3){u}\ar@{=}[d] & B\ot\syz{A}{n}N\ar[r]_(0.55){\bar{\iota}}\ar@<0.5ex>@{->>}[d] & \wbar{F}_{n-1}\ar[r]\ar@<0.5ex>@{->>}[d] & \dots \ar[r] & \wbar{F}_{1}\ar[r]\ar@<0.5ex>@{->>}[d] & \wbar{F}_{0} \ar[r] & N 
\\
N\ar[r]^(0.3){u_{1}} & B\ot\syz{A_{1}}{n{-}1}N\ar[r]^(0.6){\bar{\iota}_{1}}\ar@<0.5ex>[u] & \wbar{G}_{n-2}\ar[r]\ar@<0.5ex>[u] & \dots\ar[r] & \wbar{G}_{0}\ar[r]\ar@<0.5ex>[u] & N &
}
\end{equation}
Since \(\ob(A/J^{2}\ra B,N)=0\Ra \ob(A_{1}/(f_{2},\dots,f_{n})^{2}\ra B,N)=0\) the map \(u_{1}\) splits by induction on \(n\). So \(u\) splits. 
The other direction follows from \cite[3.6]{aus/din/sol:93}.
\end{proof}
%%%%%%%%%%%%%% PROPOSITION %%%%%%%%%%%%%%%%%
\begin{prop}\label{prop.obssplit}
Suppose \(h\co S\ra \mc{A}\) is a local Cohen-Macaulay map\textup{,} \(\mc{J}=(f_{1},\dots,f_{n})\) an \(h\)-sequence\textup{,} \(\wbar{h}\co S\ra\mc{B}=\mc{A}/\mc{J}\) the local Cohen-Macaulay map induced from \(h\)\textup{,} and \((\wbar{h},\mc{N})\) an object in \(\MCM\)\textup{.} Let 
\begin{equation*}
\xi\co\quad  0\ra\mc{L}\lra \mc{M}\xra{\;\pi\;} \mc{N}\ra 0
\end{equation*}
be the minimal \(\MCM\)-approximation of \(\mc{N}\) over \(h\)\textup{.} Then tensoring \(\xi\) by \(\mc{B}\) gives a \(4\)-term exact sequence
\begin{equation*}
0\ra \mc{N}\ot \mc{J}/\mc{J}^{2}\lra \wbar{\mc{L}}\lra \wbar{\mc{M}}\xra{\,\,\wbar{\pi}\,\,} \mc{N}\ra 0
\end{equation*}
which represents the obstruction class \(\ob(q\co\mc{A}/\mc{J}^{2}\ra\mc{B},\mc{N})\in\xt{2}{\mc{B}}{\mc{N}}{\mc{N}\ot \mc{J}/\mc{J}^{2}}\)\textup{.}

Moreover\textup{,} 
\begin{equation*}
\ob(q,\mc{N})=0 \;\Llra\; \ob(q,\mc{N}^{\vee})=0 \;\Llra\; \wbar{\pi}\tn{ splits}
\end{equation*}
where \(\mc{N}^{\vee}=\xt{n}{\mc{A}}{\mc{N}}{\omega_{h}}\)\textup{.}
\end{prop}
\begin{proof}
If $K(\ul{f})$ denotes the Koszul complex then \(\tor{\mc{A}}{i}{\mc{B}}{\mc{M}}=\cH_{i}(K(\ul{f})\ot \mc{M})=0\) for \(i>0\) by \cite[5.1-2]{ogu/ber:72}; cf. \cite[Sec. 2.2]{ile:14}. There is a map from the defining short exact sequence \(0\ra\syz{\mc{A}}{}\mc{N}\ra F_{0}\ra \mc{N}\ra 0\) to \(\xi\) extending \(\id_{\mc{N}}\). Tensoring with \(\mc{B}\) gives a map of \(4\)-term exact sequences with outer terms canonically identified. Hence they represent the same element \(\ob(q,\mc{N})\) in \(\xt{2}{\mc{B}}{\mc{N}}{\mc{N}\ot \mc{J}/\mc{J}^{2}}\).

By the argument in \cite[5.6]{ile:12} we can assume that \(\xi\) is given as \(0\ra\im(d_{n}^{\vee})\ra(\syz{\mc{A}}{n}\mc{N}^{\vee})^{\vee}\ra \mc{N}^{\vee}{}^{\vee}\ra0\) where \((F_{*},d_{*})\) is a minimal \(\mc{A}\)-free resolution of \(\mc{N}^{\vee}\). By Lemma \ref{lem.obssplit}, \(\ob(q,\mc{N}^{\vee})=0\) if and only if \(u\co \mc{N}^{\vee}\ra \mc{B}\ot\syz{\mc{A}}{n}\mc{N}^{\vee}\) splits. But applying \(\hm{}{\mc{B}}{-}{\omega_{\wbar{h}}}\) to \(u\) gives \(\wbar{\pi}\) since \(\mc{N}\cong\xt{n}{\mc{A}}{\mc{N}^{\vee}}{\omega_{h}}\cong \hm{}{\mc{B}}{\mc{N}^{\vee}}{\omega_{\wbar{h}}}\) (use 3.3.10 and the bottom of p. 114 in \cite{bru/her:98} combined with base change theory; cf. \cite[2.4]{ile:12}).
\end{proof}
%%%%%%%%%%%% REMARK %%%%%%%%%%%%
\begin{rem}\label{rem.obssplit}
In the absolute Gorenstein case with \(n=1\) this is given in \cite[4.5]{aus/din/sol:93}.
\end{rem}
\begin{proof}[Proof of Theorem {\ref{thm.defMCM}}]
Let \(\phi\) denote the composition
\begin{equation}
F:=\df{\mc{B}}{N}\lra \df{\mc{A}}{N}\xra{\,\sigma_{M}\,} \df{\mc{A}}{M}=:G
\end{equation}
Formal versality in the complete case and versality in the algebraic case implies that there is a lifting \(f\co R_{M}\ra R_{N}\) of \(\phi\).
The theorem follows from Lemma \ref{lem.res} once the conditions (i) and (ii) are verified. 

(i) By Proposition \ref{prop.obsmodule}, \(t_{F\hspace{-0.1em}/\hspace{-0.1em}\vL}\cong \xt{1}{B}{N}{N}\) and \(t_{G\hspace{-0.1em}/\hspace{-0.1em}\vL}\cong\xt{1}{A}{M}{M}\). Let \(\pi\co M\ra N\) denote the \(\MCM_{A}\)-approximation and \(\wbar{\pi}\co \wbar{M}\ra N\) the \(B\)-quotient of \(\pi\). Then \(\wbar{\pi}\) splits by Proposition \ref{prop.obssplit}. Let \(\nu\co N\ra \wbar{M}\) denote a splitting and \(\tau\co M\ra \wbar{M}\) the quotient map. Then \(\wbar{\pi}^{*}\co \xt{n}{B}{N}{N} \ra\xt{n}{B}{\wbar{M}}{N}\) splits for any \(n\). Since \(J\) is an \(M\)-regular sequence, \(\tau^{*}\co\xt{n}{B}{\wbar{M}}{N}\cong\xt{n}{A}{M}{N}\). Since \(\xt{i}{A}{M}{L}=0\) for \(i>0\), \(\pi_{*}\co \xt{n}{A}{M}{M}\cong \xt{n}{A}{M}{N}\) for \(n>0\). A diagram ensues:
\begin{equation}\label{eq.xt1}
\xymatrix@C+3pt@R-6pt@H-30pt{
& \xt{n}{A}{M}{N} & \xt{n}{A}{M}{M} \ar[l]_(0.49){\cong}^(0.47){\pi_{*}}\ar[d]^(0.47){\tau^{*}\tau_{*}} 
\\
\xt{n}{B}{N}{N} \ar@<0.5ex>[r]^{\wbar{\pi}^{*}} & \xt{n}{B}{\wbar{M}}{N} \ar[u]_{\cong}^{\tau^{*}}\ar@<0.3ex>@{->>}[l]^(0.47){\nu^{*}}\ar@<0.4ex>[r]^{\nu_{*}} & \xt{n}{B}{\wbar{M}}{\wbar{M}}\ar@<0.4ex>@{->>}[l]^(0.47){\wbar{\pi}_{*}}
}
\end{equation}
Since \(\phi_{k[\vare]}\co t_{F\hspace{-0.1em}/\hspace{-0.1em}\vL}\ra t_{G\hspace{-0.1em}/\hspace{-0.1em}\vL}\) corresponds to the composition of injective maps \((\pi_{*})^{-1}\tau^{*}\wbar{\pi}^{*}\) in \eqref{eq.xt1} for \(n=1\), \(\phi_{k[\vare]}\) is injective.

(ii) Suppose \(p\co R\ra S\) is a small surjection in \(\He\), put \(q={\id}\hot p\co\mc{A}_{R}\ra \mc{A}_{S}\), \(\wbar{q}=B\ot q\co\mc{B}_{R}\ra \mc{B}_{S}\), and so on. Suppose \(\mc{N}\) is in \(\df{\mc{B}}{N}(S)\), consider \(\mc{N}\) as \(\mc{A}_{S}\)-module and put \(\mc{M}=\sigma_{M}(\mc{N})\ra \mc{N}\); cf. \eqref{eq.ny}. There is a fixed map to \(\pi\) given in \eqref{eq.right2}. Then \(\ob(q,\mc{M})\) is contained in \(\xt{2}{A}{M}{M}\ot I\) by Lemma \ref{lem.omap} and we prove that it maps to \(\ob(\wbar{q},\mc{N})\) in \(\xt{2}{B}{N}{N}\ot I\) along the maps in \eqref{eq.xt1}.

Consider the short exact sequences of \(\mc{A}_{R}\)-modules
\begin{equation}
0 \ra  \Syz^{\mc{A}_{R}}(\mc{M}) \lra \mc{G}' \lra \mc{M} \ra 0
\end{equation}
where \(\mc{G}'\) is free. Apply \(-\ot_{R}S\) and obtain the \(4\)-term exact sequence of \(\mc{A}_{S}\)-modules
\begin{equation}\label{eq.obss}
0 \ra M\ot_{k}I \lra \Syz^{\mc{A}_{R}}(\mc{M})\ot_{R}S \lra \mc{G} \lra \mc{M} \ra 0
\end{equation}
which represents \(\ob(q,\mc{M})\in \xt{2}{\mc{A}_{S}}{\mc{M}}{M\ot I}\) by Lemma \ref{lem.obsmodule}. It splits into two short exact sequences along \(\Syz^{\mc{A}_{S}}(\mc{M})\) which is \(S\)-flat. Applying \(-\ot_{S}k\) to \eqref{eq.obss} gives a \(4\)-term exact sequence of \(A\)-modules
\begin{equation}\label{eq.obss2}
0 \ra M\ot_{k}I \lra \Syz^{\mc{A}_{R}}(\mc{M})\ot_{R}k \lra G \lra M \ra 0
\end{equation}
which represents \(\ob(q,\mc{M})\in \xt{2}{A}{M}{M\ot I}\), cf. Lemma \ref{lem.omap}. Since \(M\) is MCM and \(J\) is a regular sequence, applying \(B\ot_{A}-\) to  
\eqref{eq.obss2} gives another \(4\)-term exact sequence of \(B\)-modules
\begin{equation}\label{eq.obss3}
0 \ra \wbar{M}\ot_{k}I \lra B\ot_{A}\Syz^{\mc{A}_{R}}(\mc{M})\ot_{R}k \lra \wbar{G} \lra \wbar{M} \ra 0
\end{equation}
which represents \(\tau^{*}\tau_{*}\ob(q,\mc{M})\) in \(\xt{2}{B}{\wbar{M}}{\wbar{M}\ot I}\).
Pushout by \(\wbar{\pi}\ot\id\co \wbar{M}\ot I\ra N\ot I\) and pullback by \(\nu\co N\ra \wbar{M}\) gives the image of \(\ob(q,\mc{M})\) in \(\xt{2}{B}{N}{N\ot I}\) (cf. \eqref{eq.xt1}):
\begin{equation}\label{eq.obss4}
0 \ra N\ot_{k}I \lra E \lra Q \lra N \ra 0
\end{equation}
With a \(\mc{B}_{R}\)-free cover \(\mc{F}'\ra \mc{N}\), a similar argument gives a \(4\)-term sequence of \(B\)-modules
\begin{equation}\label{eq.obss5}
0 \ra N\ot_{k}I \lra \Syz^{\mc{B}_{R}}(\mc{N})\ot_{R}k \lra F \lra N \ra 0
\end{equation}
which represents \(\ob(\wbar{q},\mc{N})\). Since \(\wbar{M}\cong N{\oplus} X\) for some \(B\)-module \(X\), we may lift a sum of free covers to a free cover of \(\mc{M}\) and assume that \(\mc{G}'=\mc{G}'_{1}{\oplus}\mc{G}'_{2}\) with \(\mc{F}'=\mc{B}_{R}\ot\mc{G}_{1}'\). Then \(Q\cong F\oplus\Syz^{B}(X)\) and \(\Syz^{B}(\wbar{M})\cong \Syz^{B}(N)\oplus\Syz^{B}(X)\). Lifting $\wbar{\pi}$ gives a map from \eqref{eq.obss3} to \eqref{eq.obss5}. In particular there is a surjection \(B\ot_{A}\Syz^{\mc{A}_{R}}(\mc{M})\ot_{R}k\ra \Syz^{\mc{B}_{R}}(\mc{N})\ot_{R}k\) which restricts to the composition \(\wbar{M}\ot I\ra N\ot I\ra \Syz^{\mc{B}_{R}}(\mc{N})\ot_{R}k\). The induced map from $E$ to \(\Syz^{\mc{B}_{R}}(\mc{N})\ot_{R}k\) together with the projection from \(Q\) to $F$ gives a map from \eqref{eq.obss4} to \eqref{eq.obss5} which is the identity at the end terms. Thus they represent the same class in cohomology.
\end{proof}
Assume $Q=k[x_1,\dots,x_m]^{\tn{h}}$, $f\in \fr{m}_Q^2$ and put $B=Q/(f)$. Assume $N$ is a MCM $B$-module. Then there are endomorphisms $\phi$ and $\psi$ of $Q^{\oplus n}$ where $n=\dim_kN/\fr{m}N=e(B)\cdot\rk(N)$ with $\phi\psi=f{\cdot}\id=\psi\phi$ and $\coker\phi\cong N$. The pair $(\phi,\psi)$ is called a matrix factorisation of $f$ which defines $N$. Put $P=Q[t]^{\tn{h}}$, $F=f+t^2\in P$ and $A=P/(F)$. Define $G(\phi,\psi)=(\Phi,\Psi)$ where
\begin{equation}
\Phi=
\begin{pmatrix}
\phi & t
\\
-t & \psi
\end{pmatrix}
\qquad
\tn{and}
\qquad
\Psi=
\begin{pmatrix}
\psi & -t
\\
t & \phi
\end{pmatrix}
\end{equation}
are endomorphisms of $P^{\oplus 2n}$ in block-matrix notation. Then $(\Phi,\Psi)$ is a matrix factorisation of $F$ and thus defines an MCM $A$-module $\coker\Phi$ which we denote by $G(N)$. Indeed, $G$ defines a functor of stable categories $G\co\umod_B\ra\umod_A$ and was introduced by H. Kn\"orrer in \cite{kno:87}.
If $M$ a Cohen-Macaulay $A$-module of codimension $c$, put $M^{\vee}=\xt{c}{A}{M}{A}$. Note that $N^{\vee}\cong\hm{}{B}{N}{B}$ for any MCM $B$-module; cf. the bottom of p. 114 in \cite{bru/her:98}.
%%%%%%%%%%%%%%%% COROLLARY %%%%%%%%%%%%%%%%%%%%%% 
\begin{cor}\label{cor.defMCM}
For any $N$ in $\cat{MCM}_B$ there is an $\cat{MCM}_A$\textup{-}approximation
$$
0\ra A^{\oplus 2n}\lra G(N^{\vee})^{\vee}\lra N\ra 0.
$$
Put $M=G(N^{\vee})^{\vee}$\textup{.} 
Suppose $\df{\mc{B}}{N}$ and $\df{\mc{A}}{M}$ have formally versal formal families \textup{(}or versal families\textup{).} 
Then the minimal base rings \(R^{N}\) and \(R^{M}\) for $\df{\mc{B}}{N}$ and $\df{\mc{A}}{M}$ satisfy
\begin{equation*}
R^{N}\cong R^{M}\hspace{-0.2em}/J
\end{equation*}
where \(J\) is generated by elements lifting a \(k\)\textup{-}basis of the kernel of the map of dual Zariski tangent vector spaces $\phi_{k[\vare]}^*\co\xt{1}{A}{M}{M}^*\ra\xt{1}{B}{N}{N}^*$\textup{;} cf\textup{.} \textup{Lemma \ref{lem.res}.}
\end{cor}
\begin{proof}
Note first that a minimal $P$-free resolution of $N$ together with a homotopy for the multiplication with $F$ on the resolution is constructed from a minimal matrix factorisation $(\phi,\psi)$ for $N$:
\begin{equation}
\xymatrix@C+18pt@R+12pt@H-30pt{
N\ar[d]^{0}
& P^{n}\ar@{->>}[l]
\ar[dr] 
|-(0.5){
\left[
\begin{smallmatrix}
t \\ \psi
\end{smallmatrix}
\right]
}
\ar[d]^{{\cdot}F}
& P^{ n}\oplus P^{n}\ar[l]|-(0.44){
\left[
\begin{smallmatrix}
t & \phi
\end{smallmatrix}
\right]
}
\ar[dr] 
|-(0.5){
\left[
\begin{smallmatrix}
\psi & -t
\end{smallmatrix}
\right]
}
\ar[d]^{{\cdot}F}
& P^{ n}\ar[l]|-(0.46){
\left[
\begin{smallmatrix}
\phi \\
-t
\end{smallmatrix}
\right]
}
\ar[d]^{{\cdot}F}
& 0\ar[l]
\\
N
& P^{n}\ar@{->>}[l] 
& P^{ n}\oplus P^{n}\ar[l]|-(0.42){
\left[
\begin{smallmatrix}
t & \phi
\end{smallmatrix}
\right]
}
& P^{ n}\ar[l]|-(0.44){
\left[
\begin{smallmatrix}
\phi \\
-t
\end{smallmatrix}
\right]
}
& 0\ar[l]
}
\end{equation}
The Eisenbud construction \cite{eis:80} of an $A$-free resolution from these data gives:
\begin{equation}
\xymatrix@C+3pt@R+0pt@H-30pt{
N & A^{n}\ar@{->>}[l] & A^{2n}\ar[l]_(0.46){
\left[
\begin{smallmatrix}
t & \,\wbar{\phi}
\end{smallmatrix}
\right]
} & A^{2n}\ar[l]_{\wbar{\,\Phi}} & A^{2n}\ar[l]_{\wbar{\,\Psi}} & \dots\ar[l]_{\wbar{\,\Phi}}
}
\end{equation}
In particular there is a short exact sequence $0\la N\la A^n \la G(N) \la 0$. Applying $\hm{}{A}{-}{A}$ gives another short exact sequence $0\ra A^n\ra G(N)^{\vee}\ra N^{\vee}\ra 0$. This is then the $\cat{MCM}_A$-approximation of $N^{\vee}$. By local duality theory there is a canonical isomorphism $N^{\vee\vee}\cong N$; cf. \cite[3.3.10]{bru/her:98}. Thus the above construction applied to the MCM $B$-module $N^{\vee}$ gives the $\cat{MCM}_A$-approximation of $N$.

For the second part, note that $t$ is a non-zero divisor in $A$ and $A/(t)^2\cong B[t]/(t^2)$, hence $q\co A/(t)^2\thr B$ splits and $\ob{}{}(q,-)=0$. Then Theorem \ref{thm.defMCM} applies.
\end{proof}
%%%%%%%%%%%%%
\begin{ex}
Put $\wbar{M}=B\ot M$.  
By Proposition \ref{prop.obssplit}, $\ob{}{}(A/(t)^{2}\ra B,N)=0$ gives a splitting $\wbar{M}\cong N{\oplus} X$ where $X$ is stably isomorphic to $\Syz^B(N^{\vee})^{\vee}$ which (in the hypersurface case) is isomorphic to $\Syz^B(N)$. Since $\xt{1}{B}{\Syz^B(N)}{N}\cong\xt{2}{B}{N}{N}$, \eqref{eq.xt1} gives 
\begin{equation}
\xt{1}{A}{M}{M}\cong \xt{1}{B}{N}{N}\oplus\xt{2}{B}{N}{N}.
\end{equation}
Hence if $\dim_k\xt{i}{B}{N}{N}<\infty$ for $i=1,2$ then $\df{B}{N}$ and $\df{A}{M}$ have formally versal formal families for complete base rings; cf. Proposition \ref{prop.obsmodule} and \cite[2.11]{sch:68}.

If $\Spec B$ is an isolated singularity and $\Char k\neq 2$ then $\Spec A$ is an isolated singularity. Then $\df{B}{N}$ and $\df{A}{M}$ have versal elements over algebraic base rings; cf. \cite[2.4]{ess:90} and \cite[4.5]{ile:19bX}.
\end{ex}

\providecommand{\bysame}{\leavevmode\hbox to3em{\hrulefill}\thinspace}
\providecommand{\MR}{\relax\ifhmode\unskip\space\fi MR }
% \MRhref is called by the amsart/book/proc definition of \MR.
\providecommand{\MRhref}[2]{%
  \href{http://www.ams.org/mathscinet-getitem?mr=#1}{#2}
}
\providecommand{\href}[2]{#2}

%\bibliography{/Users/ile/Documents/Arbeid/Aktuelle/Generiske/Matematikk10}
\end{document}